\documentclass[11pt,a4paper, fleqn]{article}
\usepackage[utf8]{inputenc}
\usepackage[english,russian]{babel}
\usepackage{epsfig,array}

\usepackage{mathtools}
\mathtoolsset{showonlyrefs}

\usepackage{indentfirst} 
\usepackage[14pt]{extsizes}
\usepackage{amssymb}
\usepackage{amsmath}
\usepackage{amsthm}
\usepackage{pifont}
\usepackage{subfig}
\usepackage{graphicx}
\usepackage{wrapfig}
\usepackage[colorinlistoftodos]{todonotes}
\usepackage[colorlinks=true, allcolors=blue]{hyperref}
\usepackage{algorithm}
\usepackage{algpseudocode}
\usepackage{titlesec}
\titleformat{\section}
  {\large\bfseries}{\thesection.}{1em}{}
\titleformat{\subsection}
  {\centering\normalfont\normalsize\bfseries\itshape}{\thesubsection.}{1em}{}
\addto\captionsrussian{}

\newtheorem{theorem}{Теорема}

\newtheorem{Lem}{Лемма}

\newtheorem{Coro}{Утверждение}

\def\UDK#1{{\leftline{УДК {#1}}}}

\usepackage{geometry}
\def \PP {\mathbb P}
\def \EE {\mathbb E}

\newcommand{\circledOne}{\text{\ding{172}}}
\newcommand{\circledTwo}{\text{\ding{173}}}
\newcommand{\circledThree}{\text{\ding{174}}}

\def\RR{\mathbf R}
\newcommand{\argmax}{\mathop{\arg\,\max}}
\DeclarePairedDelimiter{\ceil}{\lceil}{\rceil}
\usepackage[labelsep=period]{caption}

\def\ai#1{{\color{black}#1}} 
\def\ev#1{{\color{black}#1}} 
\def\evv#1{{\color{black}#1}} 
\def\pd#1{{\color{black}#1}} 

\begin{document}
\renewcommand{\abstractname}{\vspace{-\baselineskip}}

$$
\\\\
$$

{\it \UDK{519.85}}

\begin{center}
\textbf{Численные методы для задачи распределения ресурсов в компьютерной сети}
\footnote{
Работа 
А.В.~Гасникова поддержана грантом 19-31-51001 Научное наставничество, работа П.Е.~Двуреченского поддержана грантом РФФИ 18-29-03071 мк. Работа Е.А.~Воронцовой была выполнена при поддержке Министерства науки и высшего образования Российской Федерации (госзадание) №075-00337-20-03, номер проекта 0714-2020-0005.
}

\textbf{\copyright~2021~г. 
Е.\,А.~Воронцова$^{1}$,
А.\,В.~Гасников$^{1,2}$,
П.\,Е.~Двуреченский$^{3,2}$,
А.\,С.~Иванова$^{4}$,
Д.\,А.~Пасечнюк$^{1}$
}

\textit{
$^{1}$141701 Московская обл., Долгопрудный, Институтский пер. 9, 
Московский физико-технический институт (национальный исследовательский университет), Россия \\
$^{2}$127051 Москва, Большой Каретный пер. 19, стр.1, Институт проблем передачи информации им. А.А. Харкевича РАН, Россия \\
$^{3}$Институт прикладного анализа и стохастики им. К. Вейерштрасса, Берлин, Германия \\
$^{4}$109028 Москва, Покровский бульвар 11, 
Национальный исследовательский университет <<Высшая школа экономики>>,
Россия \\
{$^*$ e-mail: asivanova@hse.ru}}

{Поступила в редакцию 29.11.2019 г.\\
Переработанный вариант  10.09.2020 г.\\Принята к публикации 00.00.2020 г.}

\end{center}

\begin{abstract}
\noindent Рассматривается задача распределения ресурсов в компьютерных сетях с большим числом соединений. Соединения
  используют для своих целей
  потребители (пользователи), число которых
  также может быть очень большим.  Для решения двойственной задачи предлагаются следующие численные методы оптимизации: быстрый градиентный метод, стохастический метод проекции субградиента, метод эллипсоидов и метод экстраполяции случайного градиента. Для каждого метода  получена оценка скорости сходимости. Также приведены алгоритмы распределенного вычисления шагов рассматриваемых методов при условии приложения их к компьютерным сетям. Отдельное внимание уделено прямо-двойственности предложенных алгоритмов.  Библ.~39. Фиг.~1. Табл.~2.

\textbf{Ключевые слова}: распределение ресурсов, сети связи, максимизация полезности сети, прямо-двойственность, быстрый градиентный метод, стохастический метод проекции субградиента, метод эллипсоидов,
метод экстраполяции случайного градиента.
\end{abstract}

\section{Введение}
\subsection{Мотивация}
Рассматривается задача управления  современными сетями связи с точки зрения оптимизации и стохастического моделирования. Для решения задач такого рода
необходимо представить
и проанализировать математическую модель, возникающую при симуляции
процесса работы крупномасштабных широкополосных сетей. Ожидается, что в будущих сетях связи появятся приложения, которые смогут изменять свои скорости передачи данных в соответствии с доступной пропускной способностью в сети. В качестве примера такой сети можно привести TCP-трафик через сеть Интернет. 

Ключевой вопрос, который мы рассматриваем в этой статье, касается того, как доступная пропускная способность в сети должна быть распределена между конкурирующими потоками. При этом 
контроль над использованием потребителями доступных пропускных способностей осуществляется посредством корректировки цены на соединение. 

Таким образом, в работе рассматривается задача
оптимизации распределения ресурсов в компьютерных сетях с большим числом соединений.
Соединения используют для своих целей
потребители (пользователи), число которых
также может быть очень большим. 
Цель работы состоит в определении механизма управления ресурсами, которые в контексте данной задачи являются доступными пропускными способностями соединений. При этом необходимо обеспечить стабильную работу системы и предотвратить перегрузки. В качестве критерия оптимальности используется сумма полезностей всех пользователей компьютерной сети.

Первоначально стандартные
задачи распределения ресурсов,
сводящиеся к максимизации
совокупной полезности производителей
при совместном использовании имеющихся
ресурсов, были рассмотрены в~\cite{kelly},
также задача распределения
ресурсов в компьютерных сетях 
исследовалась в недавней работе~\cite{rokhlin}.  Предложенные в монографии~\cite{arrow1958decentralization} механизмы децентрализованного распределения ресурсов с тех пор привлекают большое внимание в экономических исследованиях, см., например,~\cite{kakhbod2013resource}--
\cite{friedman1995complexity} и ссылки в них.  В данной работе, следуя~\cite{nesterov2018dual}--\cite{ivanova2018composite}, мы рассматриваем различные механизмы корректировки цен. Практическую ценность предложенные подходы имеют за счет своей децентрализованности, что означает,
что для установления и корректировки цены 
в отдельном соединении необходима только реакция пользователей, которые используют это соединение,
а не реакция всех пользователей сети. При
таком механизме корректировки все соединения действуют независимо. 

Кроме того, один из предложенных в статье подходов, основанный на использовании стохастического метода
проекции субградиента, обходит сле\-дую\-щую проблему, возникающую в реальных сетях. Дело в том, что в реальных сетях пакеты с данными, которые поступают от пользователей на соединение, приходят не одновременно, поэтому на практике суммарный трафик на соединении неизвестен.  
Для решения этой проблемы предлагается использовать стохастические методы, для
работы которых требуется не точное значение суммарного трафика, а только его оценка,
которую можно получить на основе трафика одного
из пользователей. \evv{Идея
использовать для решения данной проблемы стохастический метод проекции субградиента предложена
в работе~\cite{rokhlin}.}

\subsection{Содержание}
Статья организована следующим образом.
В разд.~\ref{problem_statement} описана постановка задачи и построение двойственной к ней. Также приводятся все необходимые предположения для прямой задачи. В разд.~\ref{fgm} рассматривается решение предложенной задачи быстрым градиентным методом Нестерова~\cite{Nesterov1983}, получена оценка сложности данного метода, которая имеет порядок  $O\left(\dfrac{1}{\sqrt{\varepsilon}}\right)$.  В разд.~\ref{stoch} описано решение данной задачи с использованием стохастического метода проекции субградиента и приведены оценки сложности, которые имеют порядок  $O\left(\dfrac{1}{\varepsilon^2}\right)$.  

В разд.~\ref{ellipsoids} описано решение задачи с использованием метода эллипсоидов, который хорошо подходит для задач небольшой размерности, также приведен алгоритм построения сертификата точности для данного метода. Приведены оценки сложности, которые имеют порядок  $O\left(m^2\ln\dfrac{1}{\varepsilon}\right)$,
где $m$~--- число соединений в сети. 
В разд.~\ref{reg} описана техника регуляризации, которая позволяет восстанавливать решение прямой задачи по решению двойственной для не прямо-двойственного метода. В разд.~\ref{RGEM} описано решение регуляризованной задачи с использованием метода экстраполяции случайного градиента.  Приведены оценки сложности, которые имеют порядок  $O\left(\dfrac{1}{\sqrt{\varepsilon}}\ln\dfrac{1}{\varepsilon}\right)$, где логарифмический множитель возникает за счет регуляризации двойственной задачи.

В разд.~\ref{experim} приведено 
описание вычислительных экспериментов, подтверждающих на практике теоретические результаты, полученные в предыдущих разделах. 

 Также для каждого алгоритма описано его распределенное вычисление в контексте рассматриваемой задачи.

\section{Постановка задачи}\label{problem_statement}
\begin{wrapfigure}{r}{0.5\textwidth} 
\vspace{-39pt}
  \begin{center}
    \includegraphics[width=0.5\textwidth]{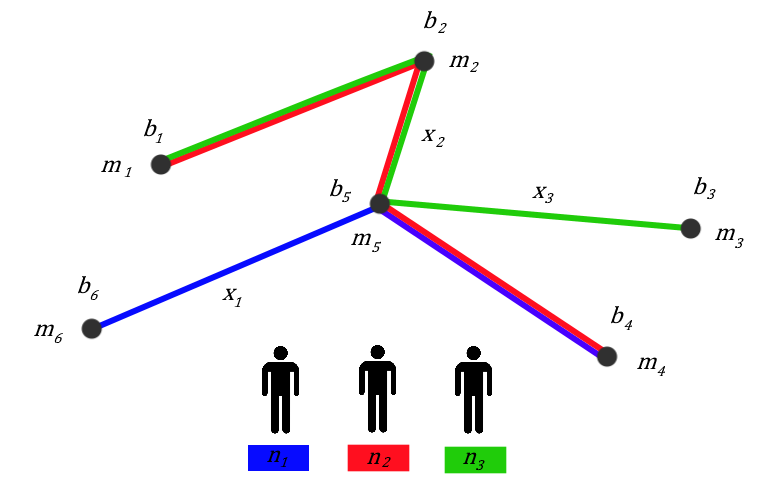}
    \caption{Пример компьютерной сети с $m=6$  и $n=3$}
    \label{fig:databaseUserTable}
  \end{center}
\vspace{-25pt}
\end{wrapfigure}

Рассмотрим компьютерную сеть с $m$ соединениями
и $n$ пользователями (или узлами),
см. фиг.~\ref{fig:databaseUserTable}. 

Пользователи
обмениваются пакетами через фиксированный набор
соединений. 
Структура сети задана матрицей маршрутизации $C = (C_i^j) \, \in \, \RR^{m \times n}$.
Столбцы матрицы $\mathbf{C}_i \neq 0$, $i = 1, \, \ldots,
\, n$~--- булевы $m$-мерные векторы, такие, что
$C_i^j = 1$ в случае использования
узлом~$i$ соединения~$j$, в противном случае $C_i^j = 0$.
Ограничения на пропускную способность
соединений задаются вектором~$\mathbf{b} \, \in \, \RR^m$
со строго положительными компонентами.

Пользователи оценивают качество работы сети
с помощью функций полезности $u_k(x_k)$, $k = 1, \, \ldots, \, n$, где $x_k \, \in \, \RR_+$~-- скорость передачи данных $k$--го пользователя.
За критерий
оптимальности системы принята сумма функций полезностей для всех пользователей~\cite{kelly}.

Задача максимизации суммарной полезности сети при заданных ограничениях на пропускную способность соединений формулируется следующим образом:
\begin{equation}
\label{eq_main_task}
    \max_{\left \{ C \mathbf{x} = \sum\limits_{k = 1}^n \mathbf{C}_k x_k \right \} \, \le \, \mathbf{b}} \left \{ U( \mathbf{x}) = \sum_{k = 1}^n u_k(x_k) \right \},
\end{equation}
где $\mathbf{x} = (x_1, \, \ldots, \, x_n) \, \in \, \mathbb{R}^n_+$.
Решением данной задачи будет оптимальное
распределение ресурсов~$\mathbf{x^*}$.

Рассмотрим стандартный переход к двойственной задаче
для~\eqref{eq_main_task}. Пусть задан вектор двойственных множителей $ \boldsymbol{\lambda} =
(\lambda_1, \, \ldots, \, \lambda_m) \, \in \, \mathbb{R}^m_+$, который можно интерпретировать как вектор цен соединений. Определим
двойственную целевую
функцию 
\begin{multline}
\varphi(\boldsymbol{\lambda}) = \max_{\mathbf{x} \, \in \, \RR_+^n} 
\left \{
\sum_{k = 1}^n u_k(x_k) + \langle \boldsymbol{\lambda},
\mathbf{b} - \sum_{k = 1}^n \mathbf{C}_k x_k \rangle
\right \} \\
= \langle \boldsymbol{\lambda}, \mathbf{b} \rangle + \sum_{k = 1}^n (u_k(x_k(\boldsymbol{\lambda})) - \langle \boldsymbol{\lambda}, \, \mathbf{C}_k x_k(\boldsymbol{\lambda}) \rangle ),
\label{eq_phi_first}
\end{multline}
причем пользователи выбирают оптимальные скорости передачи
информации $x_k$, решая следующую задачу
оптимизации
\begin{equation}
\label{eq_xi_argmax}
x_k(\boldsymbol{\lambda})  =  \argmax_{x_k \, \in \, \RR_+} \left \{
u_k(x_k) - x_k \langle \boldsymbol{\lambda}, \, \mathbf{C}_k \rangle
\right \}.
\end{equation}
Обозначим также за $\pd{\mathbf{x}}(\boldsymbol{\lambda})$ вектор с компонентами $x_k(\boldsymbol{\lambda})$.
Тогда для нахождения оптимальных цен $\boldsymbol{\lambda^*}$
требуется решить задачу
\begin{equation}
\label{dual_problem}
\min_{\boldsymbol{\lambda} \, \in \, \RR_+^m} \varphi(\boldsymbol{\lambda}).
\end{equation}

Предположим, что для прямой задачи выполняется условие Слейтера, тогда в силу сильной двойственности и прямая, и двойственная задачи будут иметь решение. Используя условие Слейтера, можно компактифицировать решение двойственной задачи. Будем предполагать, что для решения двойственной задачи верна следующая оценка: 
$$
||\boldsymbol{\lambda}^*||_2 \leq R.
$$
При этом значение $R$ никак не влияет на работу рассматриваемых алгоритмов, а только присутствует в оценках скорости их сходимости. 

Основная идея данной работы заключается в применении различных оптимизационных методов для решения двойственной задачи~\eqref{dual_problem} с добавлением прямо-двойственного анализа этих методов, позволяющего также восстановить решение прямой задачи~\eqref{eq_main_task}. В этом смысле мы развиваем подход наших предыдущих работ \cite{gasnikov2016efficient}-- 
\cite{dvurechensky2018decentralize}. Основное отличие состоит в рассмотрении ограничений-неравенств, а также в анализе стохастических алгоритмов в смысле оценок с большой вероятностью, а не в среднем. 

\subsection{Сильно вогнутые функции полезности}

В некоторых разделах мы будем предполагать, что
функции полезности $u_k(x_k)$, $k = 1, \, \ldots, \, n$, являются {\it сильно вогнутыми} с константой $\mu$. В данном подразделе опишем, какими свойствами будет обладать двойственная задача при таком предположении. 

\begin{Coro}~(Теорема Демьянова--Данскина--Рубинова,
см. \cite{Danskin}, 
\cite{Dem_book_74}).
\label{2_cor_dem}
Пусть для любого $\boldsymbol{\lambda} \, \in \, \mathbb{R}_+$ выполняется $\varphi(\boldsymbol{\lambda})=\max \limits_{\mathbf{x} \in X} F(\mathbf{x}, \, \boldsymbol{\lambda})$, где $F(\mathbf{x}, \, \boldsymbol{\lambda})$~--- выпуклая и гладкая по $\boldsymbol{\lambda}$ функция и максимум достигается в единственной точке $x(\boldsymbol{\lambda})$. Тогда $\nabla \varphi(\boldsymbol{\lambda}) = 
\evv{\nabla_{\boldsymbol{\lambda}} F(x(\boldsymbol{\lambda}), \, \boldsymbol{\lambda})}
$.
\end{Coro}

\begin{Coro} (см.~\cite{Nest}).
\label{2_cor_sopr}
Пусть для любого $k = 1, \, \ldots, \, n$  функции $u_k(x_k)$ являются $\mu$-сильно вогнутыми.
Тогда функция (\ref{eq_phi_first}), где $x_k(\boldsymbol{\lambda}),\; k =1, \, \ldots, \, n$, являются решением задачи~\eqref{eq_xi_argmax}, будет $n m^2/\mu$-гладкой, т.е. гра\-диент функции~$\varphi(\boldsymbol{\lambda})$
будет удовлетворять условию Липшица с константой $L= n m^2/\mu$: 
\begin{equation}
\label{2_smooth}
    \left|\left|\nabla \varphi(\boldsymbol{\lambda}^2)-\nabla\varphi(\boldsymbol{\lambda}^1)\right|\right|_2 \leqslant   L\left|\left|\boldsymbol{\lambda}^2 - \boldsymbol{\lambda}^1\right|\right|_2.
\end{equation}
\end{Coro}
Доказательство утверждения можно найти в Приложении.

\subsection{Вогнутые функции полезности}
Теперь предположим, что функции полезности $u_k(x_k)$, $k = 1, \, \ldots, \, n$, являются {\it вогнутыми, но не сильно вогнутыми}. Тогда двойственная задача не является гладкой. В данном подразделе описаны некоторые свойства субградиентов двойственной задачи при таких предположениях.

Субградиент двойственной задачи~\eqref{dual_problem} определяется следующим образом: 
\begin{equation*}
\ev{\nabla} \varphi(\boldsymbol{\lambda} )=
\mathbf{b} - C\mathbf{x}.
\end{equation*}
Учитывая, что $\mathbf{x}$~--- ограниченная скорость передачи информации и вектор~$\mathbf{b}$ также ограничен, получаем, что субградиенты двойственной задачи в данном случае ограничены. Таким образом, существует такая положительная константа~$M$, что верна следующая оценка:
\begin{equation}
\label{eq_gradM}
 \| \ev{\nabla}
 \varphi(\boldsymbol{\lambda} ) \|_2
\, \le \, M.
\end{equation}
\evv{В качестве грубой оценки сверху для константы~$M$ из~\eqref{eq_gradM}  можно использовать $O(n\sqrt{m})$. Множитель~$n$ возникает из-за того, что есть $n$~слагаемых, а $\sqrt{m}$ 
используется как оценка зависимости
$2$-нормы от размерности вектора~$m$. }

\section{Быстрый градиентный метод}\label{fgm}
В данном разделе мы предполагаем, что
функции полезности $u_k(x_k)$, $k = 1, \, \ldots, \, n$, являются {\it сильно вогнутыми} с константой~$\mu$, следовательно,
двойственная задача будет гладкой. 

Для решения двойственной задачи~\eqref{dual_problem} применим быстрый градиентный метод (БГМ) Нестерова 
в следующем варианте
(метод PDFGM, см. алгоритм~\ref{alg_fgm}). 
	\begin{algorithm}
		\caption{Прямо-двойственный быстрый градиентный метод (PDFGM)}\label{alg_fgm}
		\begin{algorithmic}[1]
			
			\Require $u_k(\mathbf{x}), \; \; k=1, \, \ldots, \, n$~--- сильно вогнутые функции полезности для каждого пользователя; $\boldsymbol{\lambda}^{0}$~--- вектор начальных цен, $\alpha_{t}: = \frac{t+1}{2}$, 
			\evv{$A_{-1} := 0$, } 
			$A_t := A_{t-1} + \alpha_{t} = \frac{(t+1)(t+2)}{4}$, $\tau_{t}: = \frac{\alpha_{t+1}}{A_{t+1}}=\frac{2}{t + 3}$, \evv{$t = 0, \, 1, \, \ldots, \, N-1$}.	

			\For{$t=0,\, 1, \, \ldots, \, N-1$}
			\State Вычислить $\varphi(\boldsymbol{\lambda}^t)$, $\nabla \varphi(\boldsymbol{\lambda}^t)$ 
			\State $\mathbf{y}^t :=
\left [ 
\boldsymbol{\lambda}^t - \frac{1}{L}
\left (
\mathbf{b} - \sum\limits_{k = 1}^n \mathbf{C}_k x_k (\boldsymbol{\lambda}^t)
\right )
\right ]_+$			
\State $\mathbf{z}^t := 
				\left[
				\boldsymbol{\lambda}^0 - \frac{1}{L}
		\sum\limits_{j=0}^t \alpha_j
\left (
\mathbf{b} - \sum\limits_{k = 1}^n \mathbf{C}_k x_k^{\evv{}} (\boldsymbol{\lambda}^{\evv{j}})
\right )
		\right]_{+}.
			$
			
			\State $\boldsymbol{\lambda}^{t + 1} := \tau_t \mathbf{z}^t + (1 - \tau_t) \mathbf{y}^t$
			
			\State $\mathbf{\hat{x}}^{t+1} := \frac{1}{A_{t+1}}
			\sum\limits_{j=0}^{t+1} \alpha_j \mathbf{x}(\boldsymbol{\lambda^{j}})$
\EndFor
			\State\Return $\boldsymbol{\lambda}^{N}$, $\mathbf{\hat{x}}^{N}$
		\end{algorithmic}
	\end{algorithm}

\subsection{Распределенный метод}
Для решения рассматриваемой задачи
можно также применить БГМ в распределенном варианте, что означает, что каждое соединение может вычислять оптимальную скорость передачи данных только исходя из реакции пользователей, которые используют данное соединение, и никак не взаимодействовать с другими соединениями. 

Опишем процесс, происходящий на $t$-й итерации для соединения $j$.

1. На основе информации, полученной
от пользователей
после предыдущей итерации с номером~$t-1$  (вектор $\mathbf{x}^t=\mathbf{x}(\boldsymbol{\lambda}^t)$), соединение $j$ вычисляет 
$$
y_{j}^t = \left[\lambda_{j}^t -
			\frac{1}{L}
			\left ( b_j - \sum\limits_{k = 1}^n C_k^j x_k^t \right )\right]_{+}.
			$$
При этом $C_k^j \neq 0$ только для тех пользователей, которые используют соединение~$j$. Поэтому для вычисления этого шага соединению нужна только информация от пользователей, которые используют это соединение. 

2. Далее соединение $j$ аналогично вычисляет 
$$z_{j}^t = 
			\left[\lambda_j^0 - \frac{ \alpha_j}{L}
		\left ( b_j - \sum\limits_{k = 1}^n C_k^j  {x}_k^t \right )\right]_{+}.	
			$$
			
3. Получив значения на предыдущих двух шагах, соединение~$j$ вычисляет цену 
для следующей итерации~$t+1$:
$$\lambda^{t + 1}_j = \tau_t z_j^t + (1 - \tau_t) y_j^t$$
и посылает эту информацию всем пользователям, которые с ним соединены.

4. Далее пользователи вычисляют оптимальные скорости передачи
информации $\mathbf{\hat{x}}^{t+1}$, в частности, для пользователя~$k$ получаем 
$$
x_k(\boldsymbol{\lambda}^{t+1}) = \argmax_{x_k \, \in \, \RR_+} \left (
u_k(x_k) - x_k \sum_{j = 1}^m \lambda^{t+1}_j C_k^j
\right ),
$$
где в силу определения матрицы~$C$ пользователю необходима информация только от соединений, которые он использует. Далее пользователь вычисляет оптимальную скорость
$$
\hat{x}^{t+1}_k = \frac{A_{t} \hat{x}^{t}_k + x^{t+1}_k}{A_{t+1}}.
$$

Замечание. Минусом данного алгоритма является то, что каждому соединению необходимо знать реакции всех пользователей, которые его используют на каждой итерации. К сожалению, в реальных сетях пользователи отправляют данные неодновременно, поэтому достаточно сложно собрать данную информацию для соединения. Однако наличие полной информации о пользователях позволяет соединению быстрее устанавливать равновесную цену.

\subsection{Оценка скорости сходимости БГМ}
Прежде чем привести доказательство
сходимости БГМ для рассматриваемой задачи, сформулируем ключевую лемму, необходимую для получения оценок невязки в ограничениях и зазора двойственности после
работы прямо-двойственного метода~PDFGM. 
\begin{Lem}
\label{lem_est_1}
Пусть алгоритм~$\ref{alg_fgm}$ начинает свою работу
с начальной точки $\boldsymbol{\lambda}^0$, 
которая лежит в евклидовом шаре радиуса~$R$ с центром в начале координат. Тогда после $N$ итераций алгоритма~$\ref{alg_fgm}$ будет
выполняться следующее неравенство:
\begin{equation}
\label{eq_aur37}
A_N \varphi(\mathbf{y}^N) - A_N U(\mathbf{\hat{x}}^{N})
+ 2 \hat R A_N \left \|\left ( C \mathbf{\hat{x}}^N - \mathbf{b} \right )_+ \right \|_2 \, \le \,\frac{ 37L\hat R^2}{9}.
\end{equation}
где $\mathbf{\hat{x}}^N = \frac{1}{A_N}
\sum\limits_{t = 0}^{N-1} \alpha_t \mathbf{x}(\boldsymbol{\lambda}^t)$  и $\hat R = 3R$.
\end{Lem}
Доказательство леммы можно найти в Приложении. 

Теперь сформулируем теорему об оценке скорости сходимости Алгоритма~\ref{alg_fgm}.
\begin{theorem}
\label{main_th_fgm}
Пусть Алгоритм~\ref{alg_fgm} начинает свою работу
с начальной точки $\boldsymbol{\lambda}^0$,
которая лежат в евклидовом шаре радиуса~$R$ {с центром в начале координат}.
Тогда после
\begin{equation}
\label{eq_fgm_t_iter}
N = \left \lceil \frac{2 \hat{R}}{3}\sqrt{\frac{37 L}{\varepsilon}} \right \rceil    
\end{equation}
итераций алгоритма~\ref{alg_fgm} будут
выполняться следующие неравенства:
\begin{equation*}
U(\mathbf{x^*}) - U(\mathbf{\hat{x}}^{N}) \,  \le \, \varepsilon, \, \, \, \, 
  \left \| \left (C \mathbf{\hat{x}}^N - \mathbf{b} \right)_{+}\right \|_2 \, \le \, \frac{\varepsilon}{\hat R},
\end{equation*}
где $\mathbf{\hat{x}}^N = \frac{1}{A_N}
\sum\limits_{t = 0}^{N-1} \alpha_t \mathbf{x}(\boldsymbol{\lambda}^t)$, 
$\mathbf{x^*}$ ~--- оптимальное решение задачи~\eqref{eq_main_task},
$\hat R = 3R$.
\end{theorem}
\textbf{Доказательство.}
Обозначим оптимальное значение исходной прямой задачи~\eqref{eq_main_task} $\mbox{Opt}[P]$,
а оптимальное значение двойственной задачи~\eqref{dual_problem}~--- $\mbox{Opt}[D]$.
В силу слабой двойственности
\begin{equation*}
\label{eq_dp}
\mbox{Opt}[D] \, \ge \, \mbox{Opt}[P].    
\end{equation*}
Кроме того,
для всех $\mathbf{x} \, \in \, \mathbb{R}^n_+$
и оптимального решения двойственной задачи~\eqref{dual_problem} $\boldsymbol{\lambda^*}$
\begin{equation}
\label{eq_optpu}
\mbox{Opt}[P] \, \ge \, U(\mathbf{x}) - 
\langle 
\boldsymbol{\lambda^*}, \,
\left ( \sum_{k = 1}^n \mathbf{C}_k x_k  - \mathbf{b}  \right )_+
 \rangle \geq U(\mathbf{x}) - \hat R \left \| \left ( C \mathbf{x} - \mathbf{b}   \right)_{+}\right \|_2.
\end{equation}
Тогда
\begin{eqnarray*}
\varphi(\mathbf{y}^N) -  U(\mathbf{\hat{x}}^{N}) & =& \varphi(\mathbf{y}^N) -  U(\mathbf{\hat{x}}^{N}) + \mbox{Opt}[P] - \mbox{Opt}[P] \\ &&
+ \, \mbox{Opt}[D] - \mbox{Opt}[D] \\ &
= & \underbrace{\left (\mbox{Opt}[D] - \mbox{Opt}[P] \right )}_{\geq 0}+
\left (\mbox{Opt}[P] - U(\mathbf{\hat{x}}^{N}) \right ) \\ 
&&
+ \underbrace{\left (\varphi(\mathbf{y}^N) - \mbox{Opt}[D] \right )}_{\geq 0} \, \\ &
& \overset{\eqref{eq_optpu}}{\ge} 
- \langle \boldsymbol{\lambda^*}, \, \left ( \mathbf{b} - C \mathbf{\hat{x}}^N \right )_+ \rangle 
\, \ge 
- \hat R \left \| \left ( C \mathbf{\hat{x}}^N - \mathbf{b} \right)_{+}\right \|_2.
\end{eqnarray*}
После подстановки последнего неравенства
в~\eqref{eq_aur37} получим следующую оценку: 
\begin{equation*}
 \label{eq_fgm_r37}
 \hat R \left \| \left( C \mathbf{\hat{x}}^N - \mathbf{b}\right)_{+}\right \|_2 \, \le \frac{37 L\hat R^2}{9A_N}.
\end{equation*}
Следовательно,
$\varphi(\mathbf{y}^N) -  U(\mathbf{\hat{x}}^{N}) \, \ge \, -\frac{37 L\hat R^2}{9A_N}$.
С другой стороны, из~\eqref{eq_aur37} следует,
что
$\varphi(\mathbf{y}^N) -  U(\mathbf{\hat{x}}^N) \, \le \, \frac{37 L\hat R^2}{9A_N}$.
Поэтому
$$
\left | \varphi(\mathbf{y}^N) -  U(\mathbf{\hat{x}}^{N}) \right | \, \le \, \frac{37 L\hat R^2}{9A_N}.
$$
Так как $\varphi(\mathbf{y}^N) \, \ge \, \mbox{Opt}[D] = \varphi(\mathbf{y^*}) 
\, \ge \, \mbox{Opt}[P] = U(\mathbf{x^*})$,
выполняется следующее неравенство:
\begin{equation*}
\label{eq_uu148}
U(\mathbf{x^*}) - U(\mathbf{\hat{x}}^{N}) \, \le \, \frac{37 L\hat R^2}{9A_N} = \frac{148 L\hat R^2}{9(N+1)(N+2)} \leq \frac{148 L\hat R^2}{9N^2} \leq \varepsilon.
\end{equation*}
Выражая из последнего неравенства $N$, получаем оценку из условия теоремы.

\section{Стохастический метод проекции субградиента}
\label{stoch}
Рассмотрим исходную задачу~\eqref{eq_main_task},
но теперь предположим, что
функции полезности $u_k(x_k)$, $k = 1, \, \ldots, \, n$,
{\it вогнуты, но не сильно вогнуты}.
В этом случае двойственная
задача~\eqref{dual_problem}
становится негладкой.
Поэтому для ее решения предлагается
применить стохастический метод проекции субградиента.
\evv{Впервые идея
использовать для решения данной проблемы стохастический метод проекции субградиента предложена
в работе~\cite{rokhlin}.}

Рассмотрим вероятностное пространство $\left(\Omega, \, \mathcal F, \, \PP \right)$. Пусть на нем определена последовательность независимых случайных величин $\{\xi^t \}_{t=0}^{\infty}$, равномерно распределенных на $\{1, \,
\ldots, \, n\}$, т.е.
$$
\PP(\xi^{t}=i) = \dfrac{1}{n}, \, \, \, i \, \in \, \{1, \, \ldots, \, n\}.
$$

Если доступен оракул, выдающий
стохастический субградиент двойственной функ\-ции~$\nabla \varphi(\boldsymbol{\lambda}, \, \xi)$:
$$
\nabla \varphi(\boldsymbol{\lambda}, \, \xi) =\mathbf{b}- \pd{n}C_{\xi}x_{\xi}(\boldsymbol{\lambda}), 
$$
то имеем
$$
\EE \left [\mathbf{b}- nC_{\xi_{t}}x_{\xi_{t}}(\boldsymbol{\lambda}^t) \, | \, \xi^{t} \right ] =
\mathbf{b}- \sum\limits_{k=1}^{n} \mathbf{C}_{k}x_{k}(\boldsymbol{\lambda}^t) = \nabla \varphi(\boldsymbol{\lambda}^t)
$$
Следовательно, стохастический субградиент является несмещенной оценкой субградиента.

	\begin{algorithm}
		\caption{Прямо-двойственный стохастический метод проекции субградиента~(PDSSGM), версия~1}\label{alg_stoch} 
		\begin{algorithmic}[1]
			
			\Require $u_k(\mathbf{x}), \; \; k=1, \, \ldots, \, n$~--- вогнутые функции полезности для каждого пользователя;
			$\beta$~--- шаг метода.
			\State $\boldsymbol{\lambda}^0 := \mathbf{0}$
			\For{$t= 1, \, \ldots, \, N-1$}
			\State Вычислить $\nabla \varphi(\boldsymbol{\lambda}^{t-1}, \, \xi)$
			\State $\boldsymbol{\lambda}^{t} := [\boldsymbol{\lambda}^{t-1} - \beta \left(\mathbf{b}- nC_{\xi^{t-1}}x_{\xi^{t-1}}(\boldsymbol{\lambda}^{t-1}) \right)]_+$
			\State $\hat{ \mathbf{x}}^{t+1} := \dfrac{1}{t+1}
			\sum\limits_{j=0}^{t}\mathbf{x}(\boldsymbol{\lambda}^j) $ 
			\State $\hat{\boldsymbol{\lambda}}^{t+1} := \dfrac{1}{t+1}\sum\limits_{j=0}^{t}\boldsymbol{\lambda}^j $
             \EndFor
			\State\Return $\hat{\boldsymbol{\lambda}}^{N}, \hat{ \mathbf{x}}^{N}$
		\end{algorithmic}
	\end{algorithm}

Оптимальное решение задачи~\eqref{eq_phi_first} будем
искать с помощью 
прямо-двой\-ст\-вен\-ного стохастического метода
проекции субградиента~PDSSGM. 
Приведем описание двух версий данного метода,
см. алгоритм~\ref{alg_stoch} и алгоритм~\ref{alg_stoch_1}. 
В алгоритме~\ref{alg_stoch} используется полная модель восстановления вектора~$\mathbf{x}(\boldsymbol{\lambda})$ на каждой итерации. \pd{Однако, вычисление вектора $\mathbf{x}(\boldsymbol{\lambda})$ по сложности практически эквивалентно вычислению полного субградиента $\varphi(\boldsymbol{\lambda} )$. Поэтому основным алгоритмом является  алгоритм~\ref{alg_stoch_1}, в котором} для восстановления вектора~$\mathbf{x}(\boldsymbol{\lambda})$ используется неполная, стохастическая модель, 
что означает, что на каждой итерации пересчитывается только одна компонента вектора~$\mathbf{x}(\boldsymbol{\lambda})$, остальные остаются без изменений. \pd{При доказательстве теоремы сходимости мы сначала показываем оценку скорости сходимости для алгоритма~\ref{alg_stoch}, а затем показываем, что приближенное решение прямой задачи, полученное в алгоритме~\ref{alg_stoch_1}, близко по точности к решению, полученному в алгоритме~\ref{alg_stoch}.}

\begin{algorithm}
		\caption{Прямо-двойственный стохастический метод проекции субградиента~(PDSSGM), версия~2}\label{alg_stoch_1}
		\begin{algorithmic}[1]
			
			\Require $u_k(\mathbf{x}), \; \; k=1, \, \ldots, \, n$~--- вогнутые функции полезности для каждого пользователя;
			$\beta$~--- шаг метода.
			\State $\boldsymbol{\lambda}^0 := \mathbf{0}$
			\For{$t= 1, \, \ldots, \, N-1$}
			\State Вычислить $\nabla \varphi(\boldsymbol{\lambda}^{t-1}, \, \xi)$
			\State $\boldsymbol{\lambda}^{t} := [\boldsymbol{\lambda}^{t-1} - \beta \left(\mathbf{b}- nC_{\xi^{t-1}}x_{\xi^{t-1}}(\boldsymbol{\lambda}^{t-1}) \right)]_+$
			\State $\Tilde{\mathbf{x}}^{t+1}_{\xi^t} := \frac{t}{t+1}\Tilde{\mathbf{x}}^{t}_{\xi^t} + \frac{1}{t+1}n x_{\xi^{t}}(\boldsymbol{\lambda}^{t})$, 
			\quad
			$\Tilde{\mathbf{x}}^{t+1}_{j} := \Tilde{\mathbf{x}}^{t}_j$ \mbox{ при } $j\neq \xi^t$ 
			\State $\hat{\boldsymbol{\lambda}}^{t+1} := \dfrac{1}{t+1}\sum\limits_{j=0}^{t}\boldsymbol{\lambda}^j $
             \EndFor
			\State\Return $\hat{\boldsymbol{\lambda}}^{N}, \Tilde{\mathbf{x}}^{N}$
		\end{algorithmic}
	\end{algorithm}

\subsection{Распределенный метод}
Рассмотрим, как
можно применить для решения рассматриваемой задачи
стохастический 
метод проекции субградиента в распределенном варианте. 

Опишем процесс, происходящий на $t$-й итерации для соединения~$j$:

1. На основе информации, полученной 
от соединений после предыдущей итерации с номером~$t-1$, случайный пользователь $\xi^t$ передает данные с оптимальной скоростью 
$$
x_{\xi^t}(\boldsymbol{\lambda}^{t+1}) = \argmax_{x_{\xi^t} \, \in \, \RR_+} \left (
u_{\xi^t}(x_{\xi^t}) - x_{\xi^t} \sum_{j = 1}^m \lambda^{t+1}_j C_{\xi^t}^j
\right ),
$$
где в силу определения матрицы $C$ пользователю необходима информация только от соединений, которые он использует.

2. Далее соединение $j$ вычисляет цену на следующую итерацию на основе реакции этого пользователя
$$\lambda^{t + 1}_j = \left[\lambda_{j}^t -
			\beta
			\left ( b_j - n C_{\xi^t}^j x_{\xi^t}^t \right )\right]_{+}.$$
При этом $C_{\xi^t}^j \neq 0$ только для тех пользователей, которые используют соединение $j$. Поэтому цена поменяется только для 
актуальных \evv{соединений} пользователя, который передал данные.

Замечание. { Главное преимущество данного метода заключается в том, что соединение меняет цену только на основе реакции одного пользователя, что гораздо
больше приближает
постановку задачи к ситуации
в реальных сетях, в которых пользователи отправляют данные не одновременно.}

\subsection{Оценка скорости сходимости для стохастического метода проекции субградиента}

\ai{Прежде чем перейти к доказательству основной теоремы об оценках скорости сходимости предложенных методов\ev{,} приведем необходимые предположения \ev{для решаемой задачи}.
Предположим, что существует положительная константа~$M  = O( n \sqrt{m})$, такая, что верна следующая оценка:
\begin{equation}
\label{eq_stoch_gradM}
 \| \nabla \varphi(\boldsymbol{\lambda}, \, \xi) \|_2 
\, \le \, M.
\end{equation}
Данное предположение корректно в силу того, что скорость передачи информации~$\mathbf{x}$ ограничена и вектор пропускных способностей~$\mathbf{b}$ также ограничен в силу физических соображений. Поэтому, исходя из определения стохастического субградиента, 
он ограничен. 
}

Также предположим, что
$$
\pd{\mathbb{E}}
\left [ 
\exp
\left (
\frac{\left \| \nabla \varphi(\boldsymbol{\lambda}, \, \xi) -
\nabla \varphi(\boldsymbol{\lambda})
\right \|^2_2 }{\sigma^2}
\right )
\right ]\, \le \, \exp(1),
$$
где $\sigma$~--- положительная числовая константа, \evv{порядок зависимости от $n$ и $m$ такой же, как у $M$}.

Для получения оценки скорости сходимости алгоритма~\ref{alg_stoch_1} необходимо предположение о том, что {\it функции $u_k(x_k), \, k = 1, \, \ldots, \, n$ являются липшицевыми с константой~$M_{u_k}$}, тогда функция~$U(\mathbf{x})$ будет липшицевой
с некоторой константой~$M_{U}$:
$$
\forall \, \mathbf{x}, \, \mathbf{y} \quad
   |U(\mathbf{x}) - U(\mathbf{y})| \leq M_{U}\|\mathbf{x} - \mathbf{y} \|_2,
$$
при этом $M_{U} = O(\sqrt{n})$.
Может оказаться, что функция $u_k(x_k)$ липшицева везде, кроме, например, точки~$0$. Примером такой функции является
одна из наиболее распространенных функций полезности пользователей~$u_k(x_k) = \ln x_k$. 
Но в силу специфики рассматриваемой задачи 
всегда существуют $\bar{\varepsilon} > 0$,  $\underline{\varepsilon} >0$, такие, что $x^*_k \, \geq \, \underline{\varepsilon}$, $x^*_k \, \leq \, \bar{\varepsilon}$. Тогда 
задачу можно решать на компакте~$Q = \left \{ \mathbf{x} \, : \, \underline{\varepsilon} \leq x_k \leq \bar \varepsilon, \, k = 1, \, \ldots, \, n \right \}$,
и рассматриваемая функция~$u_k(x_k) = \ln x_k$ становится липшицевой на~$Q$. В общем случае 
вогнутая функция полезности~$u(x)$ будет липшицевой
на компакте, лежащем в относительной внутренности
области определения~$u(x)$.

\ev{Пусть}
$$
\mathbb{E}
\left [ 
\exp
\left (
\frac{\left \|\mathbf{x}(\boldsymbol{\lambda}, \, \xi) -
\mathbf{x}(\boldsymbol{\lambda})
\right \|^2_2 }{\sigma_x^2}
\right )
\right ]\, \le \, \exp(1),
$$
где $\sigma_x = O(\sqrt{n})$~--- положительная числовая константа и 
$$
\mathbf{x}(\boldsymbol{\lambda}, \, \xi) = (0, \, \ldots, \, n x_{\xi}(\boldsymbol{\lambda}), \, \ldots, \, 0)^{\mbox{T}}.
$$

Сформулируем ключевую лемму, необходимую для получения оценок скорости сходимости невязки в ограничениях и зазора двойственности 
после работы прямо-двойст\-вен\-ного метода~PDSSGM.

\begin{Lem}
\label{lem_stoch_main}
Пусть алгоритм~$\ref{alg_stoch_1}$ начинает свою работу с начальной точки $\boldsymbol{\lambda}^0=\mathbf{0}$
и с шагом~$\beta$.
Тогда после $N$ итераций алгоритма~$\ref{alg_stoch_1}$
с вероятностью $1 - 4\delta$ 
выполняется неравенство:
 \begin{eqnarray*}
 \varphi(\boldsymbol{\hat{\lambda}}^N) 
&  - &  U(\Tilde{\mathbf{x}}^N)  + 
2 R
\left \|
\left [C \Tilde{\mathbf{x}}^N - \mathbf{b}
\right ]_+ 
\right \|_2   \le    C_1\dfrac{R^2 \sigma\sqrt{ g(N)J} }{\sqrt{N}} 
+ \frac{2 R^2}{\beta N} +
\frac{\beta M^2}{2}  \notag \\ & + &    \frac{ \sqrt{2} \left (1 + \sqrt{3 \ln \frac{1}{\delta}} \right )}{\sqrt{N}}  \left( M_{U} \sigma_x + 2R\left(\sigma + \sigma_x \sqrt{\lambda_{max}\left( C^{\mbox{T}}C \right)} \right) \right)\ev{,}
 \end{eqnarray*}
где 
$$
\hat{\boldsymbol{\lambda}}^N = \dfrac{1}{N} \sum\limits_{t=0}^{N-1}\boldsymbol{\lambda}^t
$$
и 
$$
\Tilde{\mathbf{x}}^N = \dfrac{1}{N} \sum\limits_{t=0}^{N-1}\mathbf{x}(\boldsymbol{\lambda}^t, \xi^{t}), 
$$ 
$C_1$-- положительная числовая константа, 
$g(N) = \ln\left(\frac{N}{\delta}\right) + \ln\ln\left(\frac{F}{f}\right)$,
$$
F  = 2\sigma^2 N (2 \beta)^{N}\left(2R^2 + 2\beta^2 M^2 +  \beta R^2 + 24 \ln\frac{N}{\delta} \beta \sigma^2 N \right),
$$
$f = \sigma^2 R^2$  и 
$$
J = \max\left\{1, \quad \frac{1}{R} \beta C_1\sqrt{\sigma^2 g(N)} + \sqrt{\frac{1}{R^2} \beta^2 C_1^2\sigma^2g(N) + \frac{2R^2 + 2\beta^2 M^2}{R^2} }\right\},
$$
\ev{а} $R$ определяется \ev{условием } $||\boldsymbol{\lambda}^*||_2 \leq R$.
\end{Lem}
Доказательство леммы приведено в Приложении.

Теперь сформулируем теорему об оценке скорости сходимости алгоритма~\ref{alg_stoch_1}.
\begin{theorem}
Пусть алгоритм~$\ref{alg_stoch_1}$ начинает свою работу
с начальной точки $\boldsymbol{\lambda}^0=\mathbf{0}$
и с шагом~$\beta =\frac{R}{M\sqrt{N}}$. Обозначим 
$$
A=  \sqrt{2} \left (1 + \sqrt{3 \ln \frac{1}{\delta}} \right )  \left( M_{U} \sigma_x + 2R\left(\sigma + \sigma_x \sqrt{\lambda_{max}\left( C^{\mbox{T}}C \right)} \right) \right) + 2.5 RM.
$$
Тогда после 
\begin{equation*}
N = O\left( \left \lceil \dfrac{A^2}{\varepsilon^2} \ln\left(\dfrac{MR}{\varepsilon\delta} \right) \right \rceil\right)
\end{equation*}
итераций алгоритма~$\ref{alg_stoch_1}$ с вероятностью $1-4\delta$ будут
выполняться следующие неравенства: 
\begin{equation*}
U(\mathbf{x^*}) - U(\Tilde{\mathbf{x}}^{N}) \,  \le \, \varepsilon, \, \, \, \, 
  \left \| \left (C \Tilde{\mathbf{x}}^N - \mathbf{b} \right)_{+}\right \|_2 \, \le \, \frac{\varepsilon}{ R},
\end{equation*}
где $\Tilde{\mathbf{x}}^N = \frac{1}{N}
\sum\limits_{t = 0}^{N-1} \mathbf{x}(\boldsymbol{\lambda}^t, \xi^{t}) $,
$\mathbf{x^*}$ ~--- оптимальное решение задачи~\eqref{eq_main_task}.
\end{theorem}
\textbf{Доказательство.}
Начало доказательства  
такое же, как в теореме~\ref{main_th_fgm}, но с использованием оценки из леммы~\ref{lem_stoch_main}. В результате для шага~$\beta=\dfrac{R}{M\sqrt{N}}$ получим, что 
\begin{eqnarray*}
 \frac{ \sqrt{2} \left (1 + \sqrt{3 \ln \frac{1}{\delta}} \right )}{\sqrt{N}}  \left( M_{U} \sigma_x + 2R\left(\sigma + \sigma_x \sqrt{\lambda_{max}\left( C^{\mbox{T}}C \right)} \right) \right)
+ 
\frac{5R M}{2\sqrt{N}}      \\
+ C_1\dfrac{
R ^ 2 \sigma\sqrt{ g(N)J} }{\sqrt{N}},
 \end{eqnarray*}
при этом с точностью до констант $g(N) \approx \ln\left(\frac{N}{\delta}\right)$, $J \approx \max
\left \{1, \, \beta\sqrt{g(N)} \right \}$. \pd{Далее найдем $N$, при котором оценка становится меньше $\varepsilon$.}

Введем следующие обозначения\ev{:}
$$
A=  \sqrt{2} \left (1 + \sqrt{3 \ln \frac{1}{\delta}} \right )  \left( M_{U} \sigma_x + 2R\left(\sigma + \sigma_x \sqrt{\lambda_{max}\left( C^{\mbox{T}}C \right)} \right) \right) + 2.5 RM, 
$$
$$
B= C_1 R^2 \sigma.
$$
Необходимо получить минимальную оценку
\ev{на} число итераций~$N$ \pd{для} достижения
заданной точности~$\varepsilon$. Для $J=1$ получаем, что 
\begin{equation}
\label{eq_stoch_teor_estN}
\sqrt{N} =\left \lceil \dfrac{A+B\sqrt{\ln\left(\frac{N}{\delta}\right)}}{\varepsilon} \right \rceil.
\end{equation}
Подставляя рекурсивно значение $N$, из
\eqref{eq_stoch_teor_estN} получим следующую оценку сложности: 
$$
N = O \left(\left \lceil \dfrac{A^2}{\varepsilon^2} \ln\left(\dfrac{MR}{\varepsilon\delta} \right) \right \rceil \right).
$$

Для $J=\beta \sqrt{g(N)}=\dfrac{R\sqrt{g(N)}}{M\sqrt{N}}$
потребуем, чтобы
$$
\dfrac{A}{\sqrt{N}} + \dfrac{B g(N) R}{M N} = \dfrac{A}{\sqrt{N}} + \dfrac{\bar B g(N) }{N} \leq \varepsilon.
$$
Так как нам нужно найти минимальное $N$, то, заменив последнее неравенство на равенство и решая полученное уравнение, получаем, что 
$$
\sqrt{N} =\left \lceil \frac{A + \sqrt{A^2 + 4\varepsilon\bar B \ln\left(\frac{N}{\delta}\right) }}{2\varepsilon} \right \rceil.
$$
Аналогично случаю $J = 1$ из последнего
равенства получается следующая оценка: 
$$
N = O \left(\left \lceil \dfrac{A^2}{\varepsilon} \ln\left(\dfrac{MR}{\varepsilon\delta} \right) \right \rceil \right).
$$
Худшая из оценок сложности для $J = 1$ и для $J=\beta \sqrt{g(N)}$ и будет оценкой из условия теоремы.

\section{Метод эллипсоидов}
\label{ellipsoids}

    В этом разделе \ev{для решения исходной задачи~\eqref{eq_main_task}
    предлагается применить метод эллипсоидов~\cite{elip_nem}}. Данный метод можно использовать \ev{при} небольшой размерности~\ev{($m$)} двойственной задачи или
    \ev{в случае, когда требуется высокая точность решения}. Данный метод является прямо-двойст\-вен\-ным, то есть по решению двойственной задачи можно восстановить решение прямой \ev{задачи}.

    Рассмотрим исходную задачу~\eqref{eq_main_task}
и двойственную к ней~\eqref{eq_phi_first}.
 Как и в предыдущем разделе, \ev{считаем}, что функции $u_k(x_k)$, $k = 1, \, 
     \ldots, \, n$, являются {\it вогнутыми, но не сильно вогнутыми}. Также предположим, что решение двойственной задачи лежит в евклидовом шаре радиуса~$R$ \ev{с центром в начале координат}, т.е. $\|\boldsymbol{\lambda}^{*}\|_2 \leq  R$.
     В качестве начальной точки метода возьмем нулевой вектор: $\boldsymbol{\lambda}^0 = 0$. Задачу будем решать на множестве:
     $$
    \Lambda_{2R} = \{\boldsymbol{\lambda} \, \in \, \mathbb{R}_{+}^m \, : \, \|\boldsymbol{\lambda}\|_2 \leq 2  R\}.
    $$ 
Привед\ev{е}м описание метода эллипсоидов (алгоритм~\ref{alg_ellipsoid}), который будем применять для решения двойственной задачи.
    
    \begin{algorithm}
		\caption{Метод эллипсоидов}\label{alg_ellipsoid}
		\begin{algorithmic}[1]
			\Require $u_k(x_k), k = 1, \, \ldots, \, n$~--- вогнутые функции полезностей
			
			\State $B_0 := 2R \cdot I_n$, \ev{$I_n$~--- единичная матрица}
			
			\For{$t=0, \, \ldots, \, N-1$}
    			\State Вычислить $\nabla \varphi(\boldsymbol{\lambda}^t)$
    			
    			\State $\mathbf{q}_t := B_t^{\mbox{T}} \nabla \varphi(\boldsymbol{\lambda}^t)$
    			\State $\mathbf{p}_t := \displaystyle \frac{B_t^{\mbox{T}} \mathbf{q}_t}{\sqrt{\mathbf{q}_t^{\mbox{T}} B_t B_t^{\mbox{T}} \mathbf{q}_t}}$
    			\State $B_{t+1} := \displaystyle \frac{m}{\sqrt{m^2 - 1}} B_t + \left(\frac{m}{m+1} - \frac{m}{\sqrt{m^2 - 1}}\right) B_t \mathbf{p}_t \mathbf{p}_t^{\mbox{T}}$
    			\State $\boldsymbol{\lambda}^{t+1} := \boldsymbol{\lambda}^t - \displaystyle \frac{1}{m+1} B_t \mathbf{p}_t$
            \EndFor
			\State\Return $\boldsymbol{\lambda}^N$
		\end{algorithmic}
	\end{algorithm}
	
Для восстановления решения прямой задачи по решению двойственной необходимо определить {\it сертификат точности}~$\xi$ для метода эллипсоидов. Напомним, что сертификатом точности называется последовательность $\xi = \{\xi_t\}_{t=0}^{N-1}$ весов, таких что
	$$
	\xi_t \geq 0, \quad \sum_{t = 0}^{N-1} \xi_t = 1.
	$$
Построение сертификата точности в нашем случае будет осуществляться в процессе работы метода эллипсоидов,
\ev{см. алгоритм~\ref{alg_cert_ellipsoid}}, общая схема которого такова~\cite{nemirovski2010accuracy}: 
\begin{enumerate}
    \item Находим <<наиболее узкую полоску>>, содержащую эллипсоид~$Q_{N}$, остающийся после итерации $N$, т.e. такой вектор $\mathbf{h}$, что 
    \ev{на $Q_{N}$ выполняется неравенство}:
\begin{equation}
\label{eq_strip}
\max_{\boldsymbol{\lambda} \in Q_{N}} \langle \mathbf{h}, \, \boldsymbol{\lambda} \rangle - \min_{\boldsymbol{\lambda} \in Q_{N}} \langle \mathbf{h}, \, \boldsymbol{\lambda} \rangle \leq 1.
\end{equation}
\ev{Для} метода эллипсоидов \ev{все $Q_N$ представлены в виде
$$
Q_N = \left \{ B_N \mathbf{z} + \boldsymbol{\lambda}^N : \mathbf{z}^{\mbox{T}} \mathbf{z} \, \le \, 1 \right \}.
$$
Тогда для решения \eqref{eq_strip} необходимо сделать} сингулярное разложение матрицы $B_{N} = U D V$, где $U$ и $V$~--- ортогональные матрицы и $D$~--- диагональная матрица
\ev{с положительными элементами на диагонали}. Далее искомый вектор~$\mathbf{h}$ определяется следующим образом: $\mathbf{h} = 1/(2 \sigma^{i_*}) \cdot U \mathbf{e}^{i_*}$, где $i_*$~--- индекс наименьшего диагонального элемента матрицы~$D$, $\sigma^{i_*}$~--- соответствующее этому индексу значение элемента матрицы $D$,  $\mathbf{e}^i$~--- векторы стандартного базиса.
    \item Для векторов $\mathbf{h}^{+} = \left [\mathbf{h}, \, -\langle \mathbf{h}, \, \boldsymbol{\lambda}^{N} \rangle \right ]$ и $\mathbf{h}^{-} = -\mathbf{h}^{+}$ находим разложения следующего вида:
$$
\mathbf{h}^{+} = \sum_{t=0}^{N-1} \nu_t [\nabla \varphi(\boldsymbol{\lambda}^t), \, -\langle\nabla \varphi(\boldsymbol{\lambda}^t), \,  \boldsymbol{\lambda}^t\rangle] + \mathbf{y}^{+}, 
$$
$$
\mathbf{h}^{-} = \sum_{t=0}^{N-1} \mu_t [\nabla \varphi(\boldsymbol{\lambda}^t), \, -\langle\nabla \varphi(\boldsymbol{\lambda}^t), \, \boldsymbol{\lambda}^t\rangle] + \mathbf{z}^{+},
$$
существование которых следует из утверждения~4.1 \cite{nemirovski2010accuracy}. Этот шаг описывают п.~6--13 рассмотренного ниже \ev{алгоритма~\ref{alg_cert_ellipsoid}}. 

    \item Из коэффициентов разложения $\nu_t$ и $\mu_t$ векторов $\mathbf{h}^{+}$ и $\mathbf{h}^{-}$, соответственно, получаем выражения для $\xi_t,\; t \in I_{N}$, где 
    $$
    I_{N} = \left\{ t \leq N-1: \, \boldsymbol{\lambda}^t \in \text{\normalfont int}\;\Lambda_{2R} \right\}.
    $$
    Коэффициенты разложения
    \ev{определяются} только для тех точек, получаемых в процессе работы \ev{алгоритма~\ref{alg_cert_ellipsoid}}, которые принадлежат допустимому множеству. 
\end{enumerate}
	
	\begin{algorithm}
		\caption{Построение сертификата точности для метода эллипсоидов}\label{alg_cert_ellipsoid}
		\begin{algorithmic}[1]
			\Require $N-1$~--- номер итерации, на которой производится вычисление сертификата точности, $\left\{B_t, \, \boldsymbol{\lambda}^t, \, \nabla \varphi(\boldsymbol{\lambda}^t)\right\}_{t=0}^{N-1}$~--- протокол работы метода эллипсоидов после $N$~итераций
			
			\If{$\nabla \varphi(\boldsymbol{\lambda}^{N-1}) = 0$}
			    \State $\xi_t := 0$ \ev{ для всех } $t = 0, \, \ldots, \, N-2$
			    \State $\xi_{N-1} := 1$
			\Else 
			\State \normalsize $\mathbf{h} := 1/(2 \sigma^{i_*}) \cdot U \mathbf{e}^{i_*}$
		    
		        \State \normalsize $\mathbf{g}_\nu := \mathbf{h}, \, \mathbf{g}_\mu := -\mathbf{h}$
		        \For{$t=0, \, \ldots, \, N-1$}
		            \State $\mathbf{q} := B_t^{\mbox{T}} \nabla \varphi(\boldsymbol{\lambda}^t)$
		            \State $\nu_t := \left[\mathbf{g}_\nu^{\mbox{T}} B_t \mathbf{q} \right]_{+} / \|\mathbf{q}\|_2^2$
		            \State $\mathbf{g}_\nu := \mathbf{g}_\nu - \nu_t \nabla \varphi(\boldsymbol{\lambda}^t)$
		            \State $\mu_t := \left[\mathbf{g}_\mu^{\mbox{T}} B_t \mathbf{q}\right]_{+} / \|\mathbf{q}\|_2^2$
		            \State $\mathbf{g}_\mu := \mathbf{g}_\mu - \mu_t \nabla \varphi(\boldsymbol{\lambda}^t)$
		        \EndFor
		        \State \normalsize $\xi_t := (\nu_t + \mu_t) \, / \sum\limits_{i \, \in \, I_{N}} (\nu_i + \mu_i), \quad t \, \in \,  I_{N}$
			\EndIf
			
			\State\Return $\left\{\xi_t\right\}_{t=0}^{N-1}$
		\end{algorithmic}
	\end{algorithm}

Замечание. Отметим, что, в отличие от БГМ и стохастического \ev{метода проекции субградиента,} в методе эллипсоидов для вычисления шагов~4--6 алгоритма~\ref{alg_ellipsoid} необходима информация о всех компонентах градиента, \ev{т.е.} нужна информация от всех пользователей. Поэтому необходим общий центр для всех соединений, который будет собирать информацию со всех соединений и выполнять эти вычисления.

Сформулируем теорему об оценк\ev{e} скорости сходимости метода эллипсоидов для решаемой задачи. 	
	\begin{theorem}{(см. \cite{nemirovski2010accuracy}).} \label{th_ellipsoid}
	Пусть алгоритм~$\ref{alg_ellipsoid}$ начинает свою работу
с начальным шаром $B_0 
= \{\boldsymbol{\lambda} \in \mathbb{R}^m: \|\boldsymbol{\lambda}\|_2 \leq 2R\}$ и сертификат точности~$\xi$ формируется в соответствии с алгоритмом~\ref{alg_cert_ellipsoid}. Тогда после 
	\begin{equation} \label{th_saga}
	    N = 2m(m+1) \ceil[\bigg]{\ln \left(\frac{32 \cdot 4 M R}{\varepsilon}\right)}
	\end{equation}
	итераций будут выполняться следующие неравенства
	$$
    U(\mathbf{x}^{*}) - U(\hat{\mathbf{x}}^N) \leq \varepsilon, \quad \left \|[C \hat{\mathbf{x}}^N - \mathbf{b}]_{+} \right \|_2 \leq \frac{\varepsilon}
    {\evv{R}},
    $$
    где $\displaystyle \hat{\mathbf{x}}^{N} = \sum_{t \in I_{N}} \xi_t \mathbf{x}(\boldsymbol{\lambda}^t), \quad I_{N} = \left\{ t \leq N-1: \, \boldsymbol{\lambda}^t \in \text{\normalfont int}\;\Lambda_{2R} \right\}$.
	\end{theorem}
Доказательство теоремы можно найти в Приложении.

\section{Регуляризация двойственной задачи}
\label{reg}
В предыдущих разделах мы рассматривали прямо-двойственные методы для решения двойственной задачи. Однако существует стандартный подход, который позволяет по решению двойственной задачи не прямо-двойственным методом воостанавливать решение прямой задачи. Ключевой идеей данного подхода является регуляризация двойственной задачи, которая гарантирует сильную выпуклос\ev{т}ь регуляризованной задачи. Далее опишем подробно предлагае\ev{м}ый подход и сформулируем леммы, связывающие решения прямой  и двойственной задачи. 

Регуляризуем функционал \eqref{eq_phi_first} по Тихонову
$$
    \varphi_\delta(\boldsymbol{\lambda}) = \varphi(\boldsymbol{\lambda}) + \frac{\delta}{2}\|\boldsymbol{\lambda}\|^2_2
$$
и вместо задачи \eqref{dual_problem} будем решать регуляризованную задачу 
\begin{equation} \label{regularized}
    \min_{\boldsymbol{\lambda} \, \in \, \RR_+^m} \varphi_{\delta}(\boldsymbol{\lambda}).
\end{equation}
Параметр $\delta$ будет оптимально подобран позже. Как и в предыдущем разделе, предполагается, что задача решается на множестве
$$
\Lambda_{2R} = \{\boldsymbol{\lambda} \in \mathbb{R}_{+}^m: \|\boldsymbol{\lambda}\|_2 \leq 2 R\}.
$$
Для полученной регуляризованной функции сформулируем следующу\ev{ю} лемму о гладкости регуляризованной задачи. 
\begin{Lem}\label{lem_reg_continuity}
Пусть функция $\varphi(\boldsymbol{\lambda})$~--- $L$-гладкая, тогда регуляризованная функция $\varphi_\delta(\boldsymbol{\lambda})$ является $(L+\delta)$-гладкой, т.е. \ev{для любых $\boldsymbol{\lambda}^1$, $\boldsymbol{\lambda}^2 \, \in \, \RR_+^m$}
\begin{equation} \label{regularized_continuity}
\left \|\nabla \varphi_\delta(\boldsymbol{\lambda}^1) - \nabla \varphi_\delta(\boldsymbol{\lambda}^2) \right \|_2 \leq (L + \delta) \left \|\boldsymbol{\lambda}^1 - \boldsymbol{\lambda}^2 \right \|_2.
\end{equation}
\end{Lem}
\textbf{Доказательство.} Градиент регуляризованной функции 
$$
\nabla \varphi_\delta(\boldsymbol{\lambda}) = \nabla \varphi(\boldsymbol{\lambda}) + \delta \boldsymbol{\lambda}.
$$
Следовательно, имеем следующую оценку:
$$
\left \|\nabla \varphi_\delta(\boldsymbol{\lambda}^1) - \nabla \varphi_\delta(\boldsymbol{\lambda}^2) \right \|_2
= \left \|\nabla \varphi(\boldsymbol{\lambda}^1) - \nabla \varphi(\boldsymbol{\lambda}^2) + \delta (\boldsymbol{\lambda}^1 - \boldsymbol{\lambda}^2) \right \|_2
$$
$$
\leq 
\left \|\nabla \varphi(\boldsymbol{\lambda}^1) - \nabla \varphi(\boldsymbol{\lambda}^2) \right \|_2 
+ \delta \left \|\boldsymbol{\lambda}^1 - \boldsymbol{\lambda}^2 \right \|_2.
$$
Откуда \ev{в силу} утверждения~\ref{2_cor_sopr} следует \eqref{regularized_continuity}. 

Также для оценки сходимости \ev{алгоритма для } прямой задачи нам понадобится следующая вспомогательная лемма
\ev{o связи} между оценкой на градиент для двойственной задачи и оценками на сходимость по функции и невязку в ограничении для прямой задачи. 
\begin{Lem}{(см.~\cite{gasnikov2016efficient}).} \label{lem_regularized_dual}
Пусть $\mathbf{x}^{*}$~--- решение прямой задачи \eqref{eq_main_task}. Тогда
\ev{выполняются следующие неравенства:}
\begin{equation} \label{eq_limits_reg1}
    \|C \mathbf{x}(\boldsymbol{\lambda}) - \mathbf{b}\|_2 \leq \|\nabla \varphi_\delta (\boldsymbol{\lambda})\|_2 + \delta \|\boldsymbol{\lambda}\|_2,
\end{equation}
\begin{equation} \label{eq_err_reg1}
     U(\mathbf{x}^{*}) - U(\mathbf{x}(\boldsymbol{\lambda})) \leq \|\nabla \varphi_\delta(\boldsymbol{\lambda})\|_2 \cdot \|\boldsymbol{\lambda}\|_2 + \delta \|\boldsymbol{\lambda}\|_2^2,
\end{equation}
где $\mathbf{x}(\boldsymbol{\lambda})$ определяется \ev{в} \eqref{eq_xi_argmax}. 
\end{Lem}
\textbf{Доказательство.}
\ev{В силу \eqref{eq_xi_argmax}}
$$
U(\mathbf{x}(\boldsymbol{\lambda})) + \langle \boldsymbol{\lambda}, \, \mathbf{b} - C \mathbf{x}(\boldsymbol{\lambda}) \rangle \geq U(\mathbf{x}^{*}) + \langle \boldsymbol{\lambda}, \, \mathbf{b} - C \mathbf{x}^{*} \rangle \geq U(\mathbf{x}^{*}).
$$
Откуда
$$
U(\mathbf{x}(\boldsymbol{\lambda})) \geq U(\mathbf{x}^{*}) - \langle \boldsymbol{\lambda}, \, \mathbf{b} - C \mathbf{x}(\boldsymbol{\lambda}) \rangle = U(\mathbf{x}^{*}) - \langle \boldsymbol{\lambda}, \, \nabla \varphi(\boldsymbol{\lambda}) \rangle.
$$
Так как $\varphi(\boldsymbol{\lambda}) = \varphi_\delta(\boldsymbol{\lambda}) - \frac{\delta}{2} \|\boldsymbol{\lambda}\|_2^2$, имеем
$$
\|\nabla \varphi(\boldsymbol{\lambda}) \|_2 = \|\nabla \varphi_\delta(\boldsymbol{\lambda}) - \delta \boldsymbol{\lambda}\|_2 \leq \|\nabla \varphi_\delta(\boldsymbol{\lambda})\| + \delta \|\boldsymbol{\lambda}\|_2.
$$
\ev{Из последнего неравенства}, учитывая соотношение $\nabla \varphi(\boldsymbol{\lambda}) = \mathbf{b} - C \mathbf{x}(\boldsymbol{\lambda})$, получаем \eqref{eq_limits_reg1}. 

Далее, оценка \eqref{eq_err_reg1} следует из 
$$
U(\mathbf{x}^{*}) - U(\mathbf{x}(\boldsymbol{\lambda})) \leq \langle \boldsymbol{\lambda}, \, \nabla \varphi(\boldsymbol{\lambda}) \rangle \leq \|\nabla \varphi(\boldsymbol{\lambda})\|_2 \cdot \|\boldsymbol{\lambda}\|_2 \leq 
$$
$$
\qquad \leq \|\boldsymbol{\lambda}\|_2 \cdot \left ( \left \|\nabla \varphi_\delta(\boldsymbol{\lambda}) \right \|_2 + \delta \left \|\boldsymbol{\lambda} \right \|_2 \right ) \leq \|\nabla \varphi_\delta(\boldsymbol{\lambda})\|_2 \cdot \|\boldsymbol{\lambda}\|_2 + \delta \|\boldsymbol{\lambda}\|_2^2. 
$$
Кроме этого, понадобится лемма о сходимости по градиенту для регуляризованной функции.
\begin{Lem} \label{lem_reqularized_ineqs}
\ev{Пусть $\boldsymbol{\lambda}_\delta^{*}$~--- решение регуляризованной двойственной задачи.}
Име\ev{e}т место следующее неравенство:
$$
   \left \|\nabla \varphi_\delta(\boldsymbol{\lambda}^{N}) \right \|_2 \leq (L + \delta) \left \|\boldsymbol{\lambda}^{N} - \boldsymbol{\lambda}_\delta^{*} \right \|_2.
$$
\end{Lem}
\textbf{Доказательство}
\ev{ немедленно следует из леммы \ref{lem_reg_continuity} и 
равенства}
$$
\nabla \varphi_\delta(\boldsymbol{\lambda}^{*}_\delta) = 0.
$$
Мы сформулировали необходимые леммы для регуляризованной задачи и в следующем разделе рассмотрим на примере применение этого подхода.

\section{Метод экстраполяции случайного градиента}
\label{RGEM}
Рассмотрим \ev{метод экстраполяции случайного градиента}~\cite{lan2018random}. 
Отметим, что данный метод не требует пересчета градиента на каждой итерации, необходимо пересчитывать только одну его компоненту на каждой итерации, что \ev{значительно} сокращает вычисления, особенно для задач большой размерности. Так как данный метод не является прямо-двойственным, то \ev{необходимо применять} алгоритм~\ref{alg_saga} к регуляризованной задаче.

    \begin{algorithm}
		\caption{Метод экстраполяции случайного градиента (RGEM)}\label{alg_saga}
		\begin{algorithmic}[1]
			\Require Параметры $\alpha$, $\eta$, $\tau$, $\{\theta_t\}_{t=1}^{N}$
			
			\State $\boldsymbol{\lambda}^0 := \mathbf{0}$
			\State $\underline{\boldsymbol{\lambda}}_i^0 := \boldsymbol{\lambda}^0$, $i = 1, \, \ldots, \, n$
			\State $y_{-1} = y_0 = \mathbf{0}$
			\For{$t=1, \, \ldots, \, N$}
    			\State 
    			\ev{Выбрать случайным образом}
    			$k_t$ из множества $\{1, \, \ldots, \, n\}$ равномерно по всем значениям
    			\State $\tilde{\mathbf{y}}^t_k := \mathbf{y}^{t-1}_k + \alpha (\mathbf{y}^{t-1}_k - \mathbf{y}^{t-2}_k), \quad k = 1, \, \ldots, \, n$
    			\State $\boldsymbol{\lambda}^t := \left[\eta \boldsymbol{\lambda}^{t-1} - \frac{1}{n} \sum\limits_{k = 1}^n \tilde{\mathbf{y}}_k^t \right]_{+} / (\delta + \eta)$\\
    			
    			\State $\underline{\boldsymbol{\lambda}}_{k_t}^t := \left (\boldsymbol{\lambda}^t + \tau \underline{\boldsymbol{\lambda}}_{k_t}^{t-1} \right ) / (1 + \tau)$
    			\State $\underline{\boldsymbol{\lambda}}_{k}^t := \underline{\boldsymbol{\lambda}}_{k}^{t-1}$, $k \in \{1, \, \ldots \, n\} \setminus \{k_t\}$\\
    			\State $\mathbf{y}^t_{k_t} := \mathbf{b} - n \mathbf{C}_{k_t} x_{k_t}(\underline{\boldsymbol{\lambda}}_{k_t}^t)$
    			\State $\mathbf{y}^t_{k} := \mathbf{y}^{t-1}_{k}$, $k \in \{1, \, \ldots, \, n\} \setminus \{k_t\}$
    			
            \EndFor
            \State $\overline{\boldsymbol{\lambda}}^N := \left(\sum\limits_{t = 0}^{N-1} \theta_t \boldsymbol{\lambda}^t\right) / \sum\limits_{t = 1}^{N} \theta_t$
			\State\Return $\overline{\boldsymbol{\lambda}}^N$
		\end{algorithmic}
	\end{algorithm}
Параметры $\alpha$, $\eta$, $\tau$ и $\theta_t$ выбираются следующим образом:
\begin{equation} \label{params1}
\overline{\alpha} = 1 - \frac{1}{n + \sqrt{n^2 + 16 n L / \delta}},
\end{equation}
\begin{equation} \label{params}
\alpha = n \overline{\alpha}, \quad \eta = \frac{\delta \overline{\alpha}}{1 - \overline{\alpha}}, \quad \tau = \frac{1}{n (1 - \overline{\alpha})} - 1, \quad
\theta_t = \overline{\alpha}^{-t}.
\end{equation}

\subsection{Распределенный метод}
В данной секции опишем распределенную версию рассмотренного метода. Для начала отметим, что векторы~$\underline{\boldsymbol{\lambda}}_1^0, \, \ldots, \, \underline{\boldsymbol{\lambda}}_n^0$ хранятся у соответствующих пользователей и влияют на формирование оптимального значения трафика для соответствующего пользователя. Как было отмечено при описании распределенного БГМ, на определение оптимального трафика пользователя влияют только цены соединений, \ev{через} которые данный пользователь \ev{обменивается пакетами}. Поэтому можно считать, что в векторе~$\underline{\boldsymbol{\lambda}}_k^t$ ненулевые только те компоненты, номера которых совпадают с номерами используемых соединений. 

Опишем распредел\ev{е}нный алгоритм на $t$-й итерации.

1. \ev{Используя} информацию, собранную с пользователей \ev{на} предыдущей итерации соединение~$j$ \ev{вычисляет} 
$$
\tilde{y}^t_{k,j} := y^{t-1}_{k,j} + \alpha 
\left ( y^{t-1}_{k,j} - y^{t-2}_{k,j} \right ) 
$$
$$
= b_{j}-nC_{k}^j x_{k} \left (\underline{\boldsymbol{\lambda}}_k^{t-1} \right ) + \alpha 
\left ( nC_{k}^j x_{k}(\underline{\boldsymbol{\lambda}}_k^{t-2}) - nC_{k}^j x_{k} \left (\underline{\boldsymbol{\lambda}}_k^{t-1} \right ) \right ).
$$
\ev{З}аметим, что в силу определения матрицы $C$ соединению $j$ нужна информация только 
\ev{от} пользователей, которые \ev{обмениваются пакетами через это соединение.} 

2. Далее соединение $j$ меняет цену по следующему правилу\ev{:} 
$$
\lambda_j^t = \left[\eta \lambda_j^{t-1} - \frac{1}{n} \sum_{k=1}^n \tilde{y}_{k, \, j}^t \right]_{+} / (\delta + \eta).
$$
    	
3. \ev{О}дин из пользователей $k_t$ реагирует на изменение цен и в качестве локального вектора цен сохраняет 
 $$\underline{\boldsymbol{\lambda}}_{k_t}^t = 
 \left (\boldsymbol{\lambda}^t + \tau \underline{\boldsymbol{\lambda}}_{k_t}^{t-1} \right) / (1 + \tau),
 $$
остальные пользователи никак не изменяют локальные цены, т.е.  $\underline{\boldsymbol{\lambda}}_{k_t}^t = \underline{\boldsymbol{\lambda}}_{k_t}^{t-1}$.

4. \ev{П}ользователь $k_t$  вычисляет 
$$
x_{k_t}^t(\underline{\lambda}_{k_t,j}^t) = \argmax_{x_k \, \in \, \RR_+} \left (
u_k(x_k) - x_k \sum_{j = 1}^m \underline{\lambda}_{k_t, \, j}^t C_k^j
\right )
$$
и отправляет эту информацию соединениям, которые использует. 

5. \ev{С}оединение $j$ обновля\ev{е}т у себя информацию для $k_t$ пользователя 
$$
y^t_{k_t,j} = b_j - n C_{k_t}^j x_{k_t}(\underline{\boldsymbol{\lambda}}_{k_t}^t).
$$
\ev{Д}анную информацию соединение будет обновлять только 
\ev{в случае, когда через него обменивается пакетами}
пользователь~$k_t$.

\subsection{Оценка скорости сходимости метода экстраполяции случайного градиента}

\ai{В данном разделе мы, \ev{так же, как в} разд.~\ref{fgm}, рассматриваем задачу~\eqref{problem_statement} с $\mu$-сильно вогнутыми функциями затрат $u_k(x_k), \, k = 1, \, \ldots, \, n$. Напомним, что в силу сильной вогнутости функций затрат двойственная задача~\eqref{dual_problem} будет гладкой с константой \ev{Липшица}~$L = \dfrac{nm^2}{\mu}$.}

Для оценки скорости сходимости \ev{метода необходима} следующая оценка для невязки по аргументу из теоремы~2.1 \cite{lan2018random}:
\begin{equation} \label{th_lan}
    \mathbb{E} \left[\|\boldsymbol{\lambda}^*_\delta - \boldsymbol{\lambda}^N\|_2^2\right] \leq \frac{4 \Delta \ev{(} \overline{\alpha} \ev{  )} ^N}{\delta},
\end{equation}
где $\Delta = \delta \|\boldsymbol{\lambda}^*_\delta - \boldsymbol{\lambda}^0\|_2^2 + \frac{B}{n \delta} + \varphi_\delta(\boldsymbol{\lambda}^0) - \varphi_\delta(\boldsymbol{\lambda}^*_\delta)$, $B = ||\mathbf{b}||_2^2$.

Используя \ev{\eqref{th_lan}}, 
можно \ev{доказать } 
теорему об оценке скорости сходимости метода непосредственно для задачи~\ev{\eqref{regularized}}. 
\begin{theorem}
    Пусть для решения регуляризованной двойственной задачи \eqref{regularized} применяется метод RGEM c параметрами \eqref{params1}, \eqref{params} и $\delta = \frac{\varepsilon}{8 R^2}$ для
    \begin{equation} \label{saga_T1}
 N = \left \lceil 2\left(n + \sqrt{n^2 + \dfrac{128nLR^2}{\varepsilon}}\right)\ln\left(\dfrac{4RA}{\varepsilon} \right) \right \rceil 
    \end{equation}
    итераций, где $A= 2\left( LR + \dfrac{\varepsilon}{8R}\right)\sqrt{6 + \dfrac{16LR^2n + 8B}{n\varepsilon}}$.
Тогда
    $$
    \mathbb{E} \left[U(\mathbf{x}^{*}) - U(\mathbf{x}(\boldsymbol{\lambda}^{N}))\right] \leq \varepsilon, \; \, \, \,  \mathbb{E} \left [ \left \| C \mathbf{x}(\boldsymbol{\lambda}^{N}) - \mathbf{b}\right \|_2\right] \leq \frac{\varepsilon}{2 R}.
    $$
\end{theorem}
\textbf{Доказательство.} 
Из леммы \ref{lem_regularized_dual} имеем оценки для невязки по ограничениям \eqref{eq_limits_reg1} и по \ev{целевой} функции \eqref{eq_err_reg1}. Ввиду предположения $\boldsymbol{\lambda} \in \Lambda_{2R}$, \ev{выполняются следующие неравенства:}
\begin{equation} \label{ineq_1}
    \|C \mathbf{x}^{N} - \mathbf{b}\|_2 \leq \|\nabla \varphi_\delta (\boldsymbol{\lambda}^{N})\|_2 + 2 \delta R,
\end{equation}
\begin{equation} \label{ineq_2}
    U(\mathbf{x}^{*}) - U(\mathbf{x}^{N}) \leq 2 R \|\nabla \varphi_\delta(\boldsymbol{\lambda}^{N})\|_2 + 4 \delta R^2,
\end{equation}
где  $\mathbf{x}^{N} = \mathbf{x}(\boldsymbol{\lambda}^N)$.
Из леммы \ref{lem_reqularized_ineqs} и неравенства \eqref{th_lan} получаем следующую оценку для $\|\nabla \varphi_\delta (\boldsymbol{\lambda}^{N})\|_2$\ev{:}
$$
    \mathbb{E} \left[ \|\nabla \varphi_\delta (\boldsymbol{\lambda}^{N})\|_2 \right] \leq 2 (L + \delta) \sqrt{\frac{\Delta}{\delta}} \ev{(}\overline{\alpha}\ev{)}^{N/2}. 
$$
Оценим $\Delta$. \ev{Для функции  $\varphi_\delta$ с липшицевым
градиентом выполняется неравенство 
$$
\varphi_\delta(\boldsymbol{\lambda}^0) - \varphi_\delta(\boldsymbol{\lambda}^*_\delta) \leq \langle \nabla \varphi_{\delta} (\boldsymbol{\lambda}^{*}_\delta), \, \boldsymbol{\lambda}^0 - \boldsymbol{\lambda}^* \rangle + \frac{L + \delta}{2} \|\boldsymbol{\lambda}^*_\delta - \boldsymbol{\lambda}^0\|_2^2.
$$
}
Учитывая\ev{,} что $\nabla \varphi_{\delta} (\boldsymbol{\lambda}^{*}_\delta)=0$, получаем  
$$
\Delta \leq  \delta \|\boldsymbol{\lambda}^*_\delta - \boldsymbol{\lambda}^0\|_2^2 + \frac{B}{n \delta} +\frac{L+\delta}{2}\|\boldsymbol{\lambda}^*_\delta - \boldsymbol{\lambda}^0\|_2^2 
\leq (6\delta + 2L)R^2 + \frac{B}{n \delta}.
$$

Пусть $\delta$ подбирается так, чтобы \ev{выполнялось} $4 \delta R^2 = \frac{\varepsilon}{2}$, \ev{тогда} $\delta = \frac{\varepsilon}{8 R^2}$. Отсюда имеем
$$
4 \delta R^2 = \frac{\varepsilon}{2}, \quad 2 \delta R = \frac{\varepsilon}{4 R}.
$$
Потребуем \ev{выполнения неравенства} $U(\mathbf{x}^{*}) - U(\mathbf{x}(\boldsymbol{\lambda}^N)) \leq \varepsilon$. \ev{Тогда, в силу~\eqref{ineq_1}, \eqref{ineq_2} должно выполняться неравенство:}
$$
\|\nabla \varphi_\delta(\boldsymbol{\lambda}^N)\|_2 \leq \frac{\varepsilon}{4 R}.
$$ 
Откуда имеем
$$
2 (L + \delta) \sqrt{\frac{\Delta}{\delta}} \ev{(} \overline{\alpha}\ev{)}^{N/2}  \leq \frac{\varepsilon}{4 R}.
$$
Учитывая, что 
$$
2 (L + \delta) \sqrt{\frac{\Delta}{\delta}}
\ev{\, \leq \, } 2\left( LR + \dfrac{\varepsilon}{8R}\right)\sqrt{6 + \dfrac{16LR^2n + 8B}{n\varepsilon}},
$$
получаем следующую оценку на число итераций: 
$$
 N = \left \lceil 2\left(n + \sqrt{n^2 + \dfrac{128nLR^2}{\varepsilon}}\right)\ln\left(\dfrac{4RA}{\varepsilon} \right) \right \rceil,
 $$
 где \ev{$A = 2\left( LR + \dfrac{\varepsilon}{8R}\right)\sqrt{6 + \dfrac{16LR^2n + 8B}{n\varepsilon}}$.}

\ai{{Замечание. Отметим, что оценку сложности алгоритма~\ref{alg_saga} можно также представить в виде $O \left ( \max\left\{n, \, \sqrt{nLR^2/\varepsilon} \right\} \cdot \ln \left( \frac{1}{\varepsilon} \right )\right)$, где логарифмический множитель появляется за сч\ev{е}т необходимости регуляриз\ev{ации} двойственной задачи. При этом на каждой итерации вычисляется только одна компонента \ev{вектора реакции пользователя на изменение цены, } 
соответственно\ev{,} арифметическая сложность операции \ev{лучше, чем при вычислении всех компонент этих векторов. }
\ev{Д}ля БГМ предположения 
\ev{для целевой функции }
аналогичны, но за сч\ev{е}т того, что на каждой итерации необходимо вычислять полный градиент\ev{,} стоимость работы алгоритма\ev{,} соответственно\ev{,} 
$O\left(n\sqrt{LR^2/\varepsilon} \right)$. Таким образом, несмотря на то\ev{,} что теоретическая оценка скорости сходимости для RGEM имеет тот же порядок\ev{,} что и для БГМ, 
на практике выигрыш получается за сч\ev{е}т более деш\ev{е}вых вычислений в рамках одной итерации.} }

\section{Вычислительные эксперименты}
\label{experim}
    Программный код для численных экспериментов был написан на языках программирования Python версии~3.6 и С++ стандарта C++14. Исходный код для экспериментов и рассмотренных в статье методов опубликован на GitHub и доступен по ссылке: \url{https://github.com/dmivilensky/network-resource-allocation}. Измерение времени работы производилось на компьютере с 2-ядерным процессором Intel Core i5-5250U с тактовой частотой 1,6 ГГц на каждое ядро, ОЗУ компьютера составляла 8 Гб.

    \subsection{Сильно выпуклые (квадратичные) функции по\-лез\-нос\-ти}
    Рассмотрим задачу \eqref{eq_main_task} для функций полезности следующего вида:
    $$
        u_k(x_k) = a_k x_k - \frac{\sigma n}{2} x_k^2, \quad a_k \sim \mathcal{U}(0, \, 100), \sigma=0.1,
    $$
    где $a_k$~--- независимые и одинаково распределённые случайные величины.
    Тогда задачу \eqref{eq_xi_argmax} можно решить явно:
    $$
    \mathbf{x}(\mathbf{\lambda}) = \frac{[\mathbf{a} - C \boldsymbol{\lambda}]_{+}}{n \sigma}.
    $$
    
    Для малого числа пользователей ($n=1500$) пропускные способности соединений выбираются одинаковыми (в данном случае $\mathbf{b} = (5, \, \ldots, \, 5)^\top$), а спрос на передачу информации равномерен ($c_{i j} = 1$ для любых $i, \, j$). Для большего числа пользователей вектор пропускных способностей генерируется случайно, так что $b_i \sim \mathcal{U}(1, 6)$. Также случайно и независимо выбираются элементы матрицы спроса, так что $c_{i j} = 1$ с вероятностью $p = 0.5$ и $c_{i j} = 0$ с вероятностью $q = 0.5$.
    
    В табл.~\ref{table_quadratic} представлены число итераций и время работы быстрого градиентного метода (FGM) и метода экстраполяции случайного градиента (RGEM) для различных конфигураций сети (с числом соединений $m$), числа пользователей $n$ и требуемой точности $\varepsilon$. В таблице выделены те случаи, в которых метод экстраполяции случайного градиента сходится к решению за меньшее время, чем быстрый градиентный метод, несмотря на большее число итераций по сравнению с БГМ. Действительно, для $n \gg 0$ число требующихся запросов оптимального решения $x_k(\boldsymbol{\lambda})$ от пользователей, в случае метода экстраполяции случайного градиента, будет меньшим, чем при использовании других алгоритмов, так как за одну итерацию метода экстраполяции случайного градиента запрос отправляется только к одному случайному пользователю.
  
    \begin{center}
    \begin{table}[ht!]
    \centering
    \begin{tabular}{|c|c|c|c|c|}
    \hline
        & \multicolumn{2}{c}{FGM} & \multicolumn{2}{|c|}{RGEM} \\ \hline
        Сеть & Итерации & Время & Итерации & Время  \\ \hline
        $m = 2$, $n=1500$, $\varepsilon = 10^{-2}$ & 350 & \textcolor{blue}{\textbf{24,5} с} & 3000 & \textcolor{blue}{\textbf{21,1} с} \\ 
        $m = 5$, $n=1500$, $\varepsilon = 10^{-2}$ & 380 & \textcolor{blue}{\textbf{42,7} с} & 6700 & \textcolor{blue}{\textbf{36,9} с} \\ 
        $m = 70$, $n=5000$, $\varepsilon = 10^{-2}$ & 400 & \textcolor{blue}{\textbf{150,0} с} & 7800 & \textcolor{blue}{\textbf{132,6} с} \\ 
        $m = 70$, $n=5000$, $\varepsilon = 10^{-3}$ & 1070 & \textcolor{blue}{\textbf{374,5} с} & 9180 & \textcolor{blue}{\textbf{283,7} с} \\ 
        $m = 100$, $n=5000$, $\varepsilon = 10^{-2}$ & 417 & \textcolor{blue}{\textbf{175,1} с} & 8200 & \textcolor{blue}{\textbf{164,0} с} \\ 
        $m = 70$, $n=7000$, $\varepsilon = 10^{-2}$ & 421 & \textcolor{blue}{\textbf{218,9} с} & 8600 & \textcolor{blue}{\textbf{206,4} с} \\ 
        $m = 100$, $n=7000$, $\varepsilon = 10^{-2}$ & 427 & \textcolor{blue}{\textbf{290,3} с} & 9200 & \textcolor{blue}{\textbf{276,0} с} \\ 
        $m = 100$, $n=7000$, $\varepsilon = 10^{-3}$ & 1120 & \textcolor{blue}{\textbf{761,6} с} & 10130 & \textcolor{blue}{\textbf{638,2} с} \\ \hline
    \end{tabular}
    \caption{Сравнение количества итераций и времени работы БГМ и метода экстраполяции случайного градиента для сильно выпуклых (квадратичных) функций полезности}
    \label{table_quadratic}
    \end{table}
    \end{center}

    \subsection{Выпуклые (логарифмические) функции полезности}
    Рассмотрим работу стохастического субградиентного метода (алгоритм \ref{alg_stoch}) и метода эллипсоидов (алгоритм \ref{alg_ellipsoid}) для функции полезности следующего вида:
    $$
        u_k(x_k) = \ln{x_k}.
    $$
    В таком случае явное решение задачи \eqref{eq_xi_argmax} выглядит следующим образом (операция $1 / \cdot$ применительно к вектору понимается поэлементно)\ev{:}
    $$
        \mathbf{x}(\boldsymbol{\lambda}) = \frac{1}{C \boldsymbol{\lambda}}.
    $$
    
    Для малого числа пользователей ($n=1500$) пропускные способности соединений выбираются одинаковыми (в данном случае $\mathbf{b} = (5, \,\
    \ldots, \,5)^\top$), а спрос на передачу информации равномерен ($c_{i j} = 1$ для любых $i, \, j$). Для большего числа пользователей вектор пропускных способностей генерируется случайно, так что $b_i \sim \mathcal{U}(1, 6)$. Также случайно и независимо выбираются элементы матрицы спроса, так что $c_{i j} = 1$ с вероятностью $p = 0.5$ и $c_{i j} = 0$ с вероятностью $q = 0.5$.
    
    В табл.~\ref{table_log} представлены число итераций и время работы стохастического субградиентного метода (SGM) и метода эллипсоидов для различных конфигураций сети, числа пользователей и требуемой точности. В таблице выделены те случаи, в которых стохастический субградиентный метод сходится к решению за меньшее время, чем метод эллипсоидов.
    
    Отметим, что аналогично
    методу экстраполяции случайного градиента, стохастический субградиентный метод требует вычисления лишь одной компоненты вектора $\mathbf{x}(\boldsymbol{\lambda})$ реакций пользователей на установившиеся цены на каждой итерации. Таким образом, при большем числе итераций метода число вычисляемых компонент $x_k(\boldsymbol{\lambda}^t)$ будет меньшим по сравнению с другими алгоритмами, например, с методом эллипсоидов, также как и коммуникационная сложность в случае распределенной реализации.
    
    \begin{center}
    \begin{table}[ht!]
    \centering
    \begin{tabular}{|c|c|c|c|c|}
    \hline
        & \multicolumn{2}{c}{Метод эллипсоидов} & \multicolumn{2}{|c|}{SGM} \\ \hline
        Сеть & Итерации & Время & Итерации & Время  \\ \hline
        $m = 2$, $n=1500$, $\varepsilon = 10^{-2}$ & 40 & \textbf{0,02} с & 2000 & \textbf{0,2} с \\ 
        $m = 5$, $n=1500$, $\varepsilon = 10^{-2}$ & 85 & \textbf{0,06} с & 2500 & \textbf{0,3} с \\ 
        $m = 70$, $n=5000$, $\varepsilon = 10^{-2}$ & 120 & \textcolor{blue}{\textbf{1,9} с} & 4000 & \textcolor{blue}{\textbf{1,3} с} \\ 
        $m = 70$, $n=5000$, $\varepsilon = 10^{-3}$ & 800 & \textcolor{blue}{\textbf{5,4} с} & 9020 & \textcolor{blue}{\textbf{2,4} с} \\ 
        $m = 100$, $n=5000$, $\varepsilon = 10^{-2}$ & 300 & \textcolor{blue}{\textbf{9,0} с} & 5000 & \textcolor{blue}{\textbf{3,1} с} \\ 
        $m = 70$, $n=7000$, $\varepsilon = 10^{-2}$ & 250 & \textcolor{blue}{\textbf{8,7} с} & 5590 & \textcolor{blue}{\textbf{5,5} с} \\ 
        $m = 100$, $n=7000$, $\varepsilon = 10^{-2}$ & 380 & \textcolor{blue}{\textbf{19,0} с} & 6480 & \textcolor{blue}{\textbf{10,8} с} \\ 
        $m = 100$, $n=7000$, $\varepsilon = 10^{-3}$ & 1830 & \textcolor{blue}{\textbf{91,5} с} & 17970 & \textcolor{blue}{\textbf{30,6} с} \\ \hline
    \end{tabular}
    \caption{Сравнение количества итераций и времени работы стохастического субградиентного метода и метода эллипсоидов для выпуклых (логарифмических) функций полезности}
    \label{table_log}
    \end{table}
    \end{center}

\clearpage

\section{Заключение}
В заключение отметим некоторые возможные развития данной работы и кратко опишем методы без подробного анализа оценки скорости сходимости. 

В разд.~\ref{ellipsoids} мы рассматривали метод эллипсоидов для задач небольшой размерности, который является прямо-двойственным.  
\ev{С}уществуют и другие методы, которые дают высокую точность и хорошо подходят для задач небольшой размерности. Одним из таких методов является метод Вайды \cite{bubeck}. Однако для восстановления решения прямой задачи при решении двойственной задачи методом Вайды, необходимо  чтобы была сходимость по норме градиента для двойственной задачи. Для этого требуется гладкость двойственной задачи, что порождается сильной выпуклостью 
\ev{целевой функции}
прямой \ev{задачи}~(утверждение~\ref{2_cor_sopr}). Если сильной выпуклости \ev{в} прямой задаче нет, то \ev{ее} можно регуляризовать,  
как  
было описано в разд.~\ref{reg}, но при этом оценка сходимости логарифмически ухудшится. 

Также для решения двойственной задачи можно применить методы высокого порядка \pd{\cite{nesterov2018implementable}--\cite{gasnikov2019near}}, если двойственная задача достаточно гладкая. \ai{При этом шаги данных методов можно будет считать распределенно}, так как в данной задаче рассматривается централизованная архитектура с точки зрения взаимодействия соединения и использующих его пользователей.   
Однако отметим, что оптимальные методы высокого порядка, требующие одномерного поиска и не обладающие прямо-двойственностью, необходимо применять \ev{только после } 
регуляризации двойственной задачи. 

Ещ\ev{е} одно направление\ev{~---} методы типа редукции дисперсии, например, \cite{zhou2018simple}--\cite{zhou2018direct} которые являются промежуточными между \ev{стохастическим градиентным методом} и \ev{БГМ}. Однако данные методы \ev{тоже} не прямо-двойст\-вен\-ные, поэтому их необходимо применять к предварительно регуляризованной двойст\-вен\-ной задаче. 

\ev{О}собый интерес для авторов представляют метод Hogwild!~\cite{recht2011hogwild} и методы с использованием мини-батчинга. Заметим, что
\ev{при этом} одновременно данные передают не все пользователи, но и не только один пользователь, как это получается при использовани\ev{и} стохастических методов. \ev{В}ыбирая в качестве размера батча количество пользователей, которые отправляют данные в один момент времени, можно учесть специфику работы реальных сетей.



\appendix
\renewcommand{\theequation}{П.\arabic{equation}}
\section*{{$\qquad \qquad \qquad \qquad \qquad \qquad \qquad \qquad \qquad$\itПРИЛОЖЕНИЕ} } \label{sec:apCATD}

\subsection*{Вспомогательные результаты}

Привед\ev{е}м некоторые леммы из других работ, на которые мы будем ссылаться в доказательствах. Также привед\ev{е}м доказательства утверждений о свойствах двойственной функции, которые используются \ev{при} доказательстве основных теорем. 

\begin{Lem} (см. \cite[лемма~2]{jin2019short}).\label{lem:jin_lemma_2}
    Для случайного вектора  $\xi \, \in  \, \RR^n$ следующие утверждения эквивалентны с точностью до константного множителя у  $\sigma$.
    \begin{enumerate}
        \item Хвосты: $\PP\left\{\|\xi\|_2 \ge \gamma\right\} \le 2 \exp\left(-\frac{\gamma^2}{2\sigma^2}\right)$ $\forall \gamma \ge 0$.
        \item Моменты: $\left(\EE\left[\xi^p\right]\right)^{\frac{1}{p}} \le \sigma\sqrt{p}$ для любого положительного целого $p$.
        \item  Предположение легких хвостов: $\EE\left[\exp\left(\frac{\|\xi\|_2^2}{\sigma^2}\right)\right] \le \exp(1)$.
    \end{enumerate}
\end{Lem}

\begin{Lem}
(см. \cite[следствие~8]{jin2019short}).
\label{lem:jin_corollary}
    Пусть $\{\xi^k\}_{k = 1}^N$~\ev{---} последовательность случайных векторов из $\RR^n$\ev{,} таки\ev{x,} что для  $k=1, \, \ldots, \, N$ и для любых $\gamma \ge 0$
    $$
        \EE\left[\xi^k\mid \xi^1, \, \ldots, \, \xi^{k-1}\right] = 0,\quad \EE\left[\|\xi^k\|_2 \ge \gamma \mid \xi^1, \, \ldots, \, \xi^{k-1}\right] 
        $$
        $$
        \le \, \exp\left(-\frac{\gamma^2}{2\sigma_k^2}\right)
        \quad \text{почти наверное,}
    $$
    где $\sigma_k^2$ принадлежит $\sigma(\xi^1, \, \ldots, \, \xi^{k-1})$ для всех $k=1, \, \ldots, \, N$. Пусть $S_N = \sum\limits_{k=1}^N\xi^k$. Тогда существует константа $C_1$\ev{,} такая что для любых фиксированных $\delta > 0$ и $B > b > 0$ с вероятностью $1 - \delta$ \ev{выполняется}:
    \begin{equation*}
        \text{либо } \sum\limits_{k=1}^N\sigma_k^2 \ge B\ev{,} 
        \end{equation*}
        \begin{equation*}
        \text{либо} \quad \|S_N\|_2 \le C_1\sqrt{\max\left\{\sum\limits_{k=1}^N\sigma_k^2, \, b\right\}\left(\ln\frac{2n}{\delta} + \ln\ln\frac{B}{b}\right)}.
    \end{equation*}
\end{Lem}

\begin{Lem}
(см. \cite[следствие из теор.~2.1, случай  (ii)]{JudNem}).
    Пусть $\{\xi^k\}_{k = 1}^N$~\ev{---} последовательность случайных векторов из $\RR^n$ \ev{удовлетворяет условию}
    \begin{equation*}
        \EE\left[\xi^k\mid \xi^1, \, \ldots, \, \xi^{k-1}\right] = 0 \text{ почти наверное,}\quad k=1, \, \ldots, \, N
    \end{equation*}
    и пусть $S_N = \sum\limits_{k=1}^N\xi^k$. Пусть последовательность $\{\xi^k\}_{k = 1}^N$ удовлетворяет ``light-tail''\ev{-}условию:
    \begin{equation*}
        \EE\left[\exp\left(\frac{\|\xi^k\|_2^2}{\sigma_k^2}\right)\mid \xi^1, \, \ldots, \, \xi^{k-1}\right] \, \le \, \exp(1) \text{ почти наверное,}\quad k = 1, \, \ldots, \, N,
    \end{equation*}
    где $\sigma_1,\ldots,\sigma_N$~\ev{---} положительные числа. Тогда для всех $\gamma \, \ge \, 0$ \ev{выполняется}
    \begin{equation}
        \PP\left\{\|S_N\| \ge \left(\sqrt{2} + \sqrt{2\gamma}\right)\sqrt{\sum\limits_{k=1}^N\sigma_k^2}\right\} \le \exp\left(-\frac{\gamma^2}{3}\right).
    \end{equation}
\end{Lem}

\textbf{Доказательство утверждения \ref{2_cor_sopr}}

Представим двойственную функцию в следующем виде

\begin{eqnarray*}
    \varphi(\boldsymbol{\lambda})= \sum_{k = 1}^{n}  \left \{ u_k(x_k(\boldsymbol{\lambda})) - \langle\boldsymbol{\lambda}, \, \mathbf{C}_k \rangle x_k(\boldsymbol{\lambda}) + \dfrac{1}{n}\langle \boldsymbol{\lambda}, \, \mathbf{b} \rangle \right \} =\sum_{k = 1}^{n} \varphi_k(\boldsymbol{\lambda}),
\end{eqnarray*}
при этом из утверждения \ref{2_cor_dem} получаем, что
\begin{eqnarray*}
    \nabla \varphi(\boldsymbol{\lambda}) = \sum_{k = 1}^{n} \nabla \varphi_k(\boldsymbol{\lambda}) = \sum_{k = 1}^{n} \left( \frac{1}{n} \mathbf{b} -\mathbf{C}_k x_k(\boldsymbol{\lambda})  \right).
\end{eqnarray*}
Определим 
$$
    x_k(\boldsymbol{\lambda}\ev{^1})  =  \argmax_{x_k \, \in \, \RR_+} \left \{
    u_k(x_k) - x_k \langle \boldsymbol{\lambda}\ev{^1}, \mathbf{C}_k \rangle
    \right \}, \quad x_k(\boldsymbol{\lambda}\ev{^2}) 
    $$
    $$
    =  \argmax_{x_k \, \in \, \RR_+} \left \{
    u_k(x_k) - x_k \langle \boldsymbol{\lambda}\ev{^2}, \mathbf{C}_k \rangle
    \right \}.
$$
Запишем необходимые условия максимума первого порядка \begin{eqnarray*}
    & \left \langle \nabla u_k(x_k(\boldsymbol{\lambda}\ev{^1})) - \langle \boldsymbol{\lambda}\ev{^1}, \, \mathbf{C}_k \rangle, \, x_k(\boldsymbol{\lambda}\ev{^1}) - x_k(\boldsymbol{\lambda}\ev{^2}) \right \rangle \geq 0 
    , \\
    &\left \langle \nabla u_k(x_k(\boldsymbol{\lambda}\ev{^2})) - \langle \boldsymbol{\lambda}\ev{^2}, \, \mathbf{C}_k \rangle, \, x_k(\boldsymbol{\lambda}\ev{^2}) - x_k(\boldsymbol{\lambda}\ev{^1}) \right \rangle \geq 0.
\end{eqnarray*}
Складывая эти неравенства, получаем 
$$
\left    \langle \nabla u_k(x_k(\boldsymbol{\lambda}\ev{^2})) - \nabla u_k(x_k(\boldsymbol{\lambda}\ev{^1})), \, x_k(\boldsymbol{\lambda}\ev{^1}) -  x_k(\boldsymbol{\lambda}\ev{^2}) \right \rangle 
$$
$$
\leq 
    \left   \langle \langle \boldsymbol{\lambda}\ev{^2}, \, \mathbf{C}_k \rangle  - \langle \boldsymbol{\lambda}\ev{^1}, \, \mathbf{C}_k \rangle,  \, x_k(\boldsymbol{\lambda}\ev{^1}) -  x_k (\boldsymbol{\lambda}\ev{^2}) \right \rangle.
$$
\ev{В силу} сильной вогнутости $u_k(x_k)$ \ev{для любых $x_k^1$ и
$x_k^2$, $k = 1, \, \ldots, \, n$, выполняется}
\begin{eqnarray*}
    \langle \nabla u_k(x_k^2) - \nabla u_k(x_k^1), \, x_k^1 -  x_k^2 \rangle \geq \mu ||x_k^1 - x_k^2||^2_2.
\end{eqnarray*}
Отсюда получаем, что 
\begin{eqnarray*}
    \mu ||x_k(\boldsymbol{\lambda}\ev{^1}) - x_k(\boldsymbol{\lambda}\ev{^2})||^2_2 &\leq &
       \langle \langle \boldsymbol{\lambda}\ev{^2}, \, \mathbf{C}_k \rangle  - \langle \boldsymbol{\lambda}\ev{^1}, \, \mathbf{C}_k \rangle,  \, x_k(\boldsymbol{\lambda}\ev{^1}) -  x_k(\boldsymbol{\lambda}\ev{^2}) \rangle \\ &\leq& ||\mathbf{C}_k ||_2 \cdot ||\boldsymbol{\lambda}\ev{^1} - \boldsymbol{\lambda}\ev{^2}||_2 \cdot ||x_k(\boldsymbol{\lambda}\ev{^1}) - x_k(\boldsymbol{\lambda}\ev{^2})||_2.
\end{eqnarray*}
Тогда можно получить следующую оценку \ev{для всех $\nabla \varphi_k$} градиента 
$$
    ||\nabla \varphi_k(\boldsymbol{\lambda}\ev{^1})-\nabla \varphi_k(\boldsymbol{\lambda}\ev{^2})||_2 \leq ||\mathbf{C}_k ||_2 \cdot ||x_k(\boldsymbol{\lambda}\ev{^1}) - x_k(\boldsymbol{\lambda}\ev{^2})||_2 
    $$
    $$
    \leq \frac{1}{\mu}
 ||\mathbf{C}_k ||^2_2 \cdot ||\boldsymbol{\lambda}\ev{^1} - \boldsymbol{\lambda}\ev{^2}||_2.
$$
\ev{Для матрицы $C$ c учетом ее структуры верна оценка }
$||\mathbf{C}_k ||_2 \leq m$. Тогда 
для градиента двойственной функции
\begin{eqnarray*}
    ||\nabla \varphi(\boldsymbol{\lambda}\ev{^1})-\nabla \varphi(\boldsymbol{\lambda}\ev{^2})||_2 \leq \sum_{k = 1}^{n}
    ||\nabla \varphi_k(\boldsymbol{\lambda}\ev{^1})-\nabla \varphi_k(\boldsymbol{\lambda}\ev{^2})||_2 \leq \frac{m^2n}{\mu} ||\boldsymbol{\lambda}\ev{^1} - \boldsymbol{\lambda}\ev{^2}||_2.
\end{eqnarray*}

\subsection*{Доказательство леммы~\ref{lem_est_1}}

Для доказательства леммы~\ref{lem_est_1} сначала сформулируем и докажем одну техническую лемму. 

Обозначим $\evv{d_L}
(\boldsymbol{\lambda}) = 
\frac{L}{2} \| \boldsymbol{\lambda} - \boldsymbol{\lambda^0} \|_2^2$ и
рассмотрим последовательности
$$
l_t(\boldsymbol{\lambda}) = \sum\limits_{j = 0}^t \alpha_j \left [ 
\varphi(\boldsymbol{\lambda}^j) + \langle
			\nabla \varphi (\boldsymbol{\lambda}^j), \boldsymbol{\lambda} - \boldsymbol{\lambda}^j\rangle 
			\right ] 
$$
и
$$
\psi_{\evv{t}
}(\boldsymbol{\lambda}) = l_t(\boldsymbol{\lambda}) + \evv{d_L}
(\boldsymbol{\lambda}), \evv{t = 0, \, 1, \, \ldots,}
$$
где $\{ \boldsymbol{\lambda}^j \}_{j \, \ge \, 0}$~--- последовательность точек, генерируемых алгоритмом~\ref{alg_fgm}.

\begin{Lem}\label{lem_fgm_pd}
После $N$ шагов алгоритма~\ref{alg_fgm}
выполняется следующее неравенство:
\begin{equation}
\label{eq_lem1}
A_N \varphi(\mathbf{y}^N) \, \le \, \min_{\boldsymbol{\lambda} \, \in \, \RR_{+}^m}
\psi_N (\boldsymbol{\lambda}) = \psi_N(\mathbf{z}^N).
\end{equation}
\end{Lem}

\textbf{Доказательство леммы \ref{lem_fgm_pd}}

Докажем по индукции, что \eqref{eq_lem1} верно. 
При $t = 0$ неравенство~\eqref{eq_lem1} выполняется.
\evv{Действительно, 
\begin{eqnarray}
\psi_0 & = & \min_{\boldsymbol{\lambda} \, \in \, \RR_{+}^m}
\left \{
\alpha_0 \left [ 
\varphi(\boldsymbol{\lambda}^0) + \langle
			\nabla \varphi (\boldsymbol{\lambda}^0), \boldsymbol{\lambda} - \boldsymbol{\lambda}^0\rangle 
			\right ] +
			\frac{L}{2} \| \boldsymbol{\lambda} - \boldsymbol{\lambda^0} \|_2^2
\right \}
 \notag \\
& \overset{\circledOne}{\ge} &
\alpha_0
\min_{\boldsymbol{\lambda} \, \in \, \RR_{+}^m}
\left \{ 
\varphi(\boldsymbol{\lambda}^0) + \langle
			\nabla \varphi (\boldsymbol{\lambda}^0), \boldsymbol{\lambda} - \boldsymbol{\lambda}^0\rangle 
			 +
			\frac{L}{2} \| \boldsymbol{\lambda} - \boldsymbol{\lambda^0} \|_2^2
\right \} 
\, \overset{\circledTwo}{\ge} \, \alpha_0 \varphi(y_0),
\nonumber
\end{eqnarray}
где
$\circledOne$ выполняется, так как 
$\alpha_0 = 1/2 \, \le \, 1$,
а $\circledTwo$~--- в силу того, что функция
$\varphi(\boldsymbol{\lambda})$ имеет липшицев
градиент (см. утверждение~\ref{2_cor_sopr}
и \cite[лемма 1.2.3]{nesterov_book2018}).
Итак, $A_0 \varphi(\mathbf{y}^0) = \frac{1}{2}
\varphi(\mathbf{y}^0) \, \le \, \psi_0$.
}

Пусть \eqref{eq_lem1} верно при $t$:
\begin{equation}
\label{eq_ind_pred}
A_t \varphi(\mathbf{y}^t) \, \le \, \psi_t(\mathbf{z}^t).
\end{equation}
Докажем, что \eqref{eq_lem1} верно при $t+1$.
Действительно, имеем,
\begin{eqnarray}
\label{eq_last_in_lem1}
\psi_{t+1}(\mathbf{z}^{t+1}) & =& \min_{\boldsymbol{\lambda} \, \in \, \RR_{+}^m}
\left \{
\psi_t(\boldsymbol{\lambda}) + \alpha_{t+1}
\left [ 
\varphi(\boldsymbol{\lambda}^{t+1}) + 
			\langle	\nabla \varphi (\boldsymbol{\lambda}^{t+1}), \boldsymbol{\lambda} - \boldsymbol{\lambda}^{t+1}\rangle 
			\right ] 
\right \} \notag \\
&
\overset{\circledOne}{\ge}& \, 
\min_{\boldsymbol{\lambda} \, \in \, \RR_{+}^m}
\Bigg \{
\psi_t(\mathbf{z}^t) + \frac{L}{2} \| \boldsymbol{\lambda} - \mathbf{z}^t \|_2^2 
\notag \\
& + & \alpha_{t+1}
\left [ 
\varphi(\boldsymbol{\lambda}^{t+1}) + 
		\langle	\nabla \varphi (\boldsymbol{\lambda}^{t+1}), \boldsymbol{\lambda} - \boldsymbol{\lambda}^{t+1}\rangle 
			\right ] 
\Bigg \} \notag  \\
&
\overset{\circledTwo}{\ge}& \,
\min_{\boldsymbol{\lambda} \, \in \, \RR_{+}^m}
\Bigg \{
A_t \varphi(\mathbf{y}^t) + \frac{L}{2} \| \boldsymbol{\lambda} - \mathbf{z}^t \|_2^2 \notag \\
& + &
\alpha_{t+1}
\left [ 
\varphi(\boldsymbol{\lambda}^{t+1}) + 
		\langle	\nabla \varphi (\boldsymbol{\lambda}^{t+1}), \boldsymbol{\lambda} - \boldsymbol{\lambda}^{t+1}\rangle 
			\right ] 
\Bigg \} \notag 
\\ &
\overset{\circledThree}{\ge} &\,
\min_{\boldsymbol{\lambda} \, \in \, \RR_{+}^m}
\Bigg \{
A_t \left ( \varphi(\boldsymbol{\lambda}^{t+1}) + \langle	\nabla \varphi (\boldsymbol{\lambda}^{t+1}), \boldsymbol{\lambda} - \boldsymbol{\lambda}^{t+1}\rangle  \right ) +
\frac{L}{2} \| \boldsymbol{\lambda} - \mathbf{z}^t \|_2^2 \notag  \\ & +&
\alpha_{t+1}
\left [ 
\varphi(\boldsymbol{\lambda}^{t+1}) 
 + \langle	\nabla \varphi (\boldsymbol{\lambda}^{t+1}), \boldsymbol{\lambda} - \boldsymbol{\lambda}^{t+1}\rangle 
			\right ] 
\Bigg \},
\end{eqnarray}
где $\circledOne$ выполняется в силу \evv{сильной выпуклости прокс-функции
$\frac{1}{2} \| \boldsymbol{\lambda} - \boldsymbol{\lambda^0} \|_2^2$
и свойств экстремума в точке $\mathbf{z}^t$},
$\circledTwo$ следует из~\eqref{eq_ind_pred},
$\circledThree$~--- в силу выпуклости функции $\varphi(\boldsymbol{\lambda})$.

Так как между коэффициентами $A_t$ и $\alpha_t$ БГМ есть следующая зависимость:
$A_{\evv{t+1}} = \sum\limits_{j=0}^{\evv{t+1}} \alpha_{\evv{j}} = A_t + \alpha_{t+1}$ и $\tau_t = \alpha_{t+1} / A_{t+1}$, соотношение $\boldsymbol{\lambda}^{t + 1} = \tau_t \mathbf{z}^t + (1 - \tau_t) \mathbf{y}^t$ из алгоритма~\ref{alg_fgm} можно переписать
в виде:
$$
A_{t+1} \boldsymbol{\lambda}^{t+1} = \alpha_{t+1} \mathbf{z}^t + A_t \mathbf{y}^t.
$$
Используя последние соотношения, можно сделать следующие
преобразования:
\begin{eqnarray*}
& A_t  \langle  \nabla \varphi (\boldsymbol{\lambda}^{t+1}), \mathbf{y}^t - \boldsymbol{\lambda}^{t+1} \rangle  + \alpha_{t+1} \langle \nabla \varphi(\boldsymbol{\lambda}^{t+1}) , \boldsymbol{\lambda} - \boldsymbol{\lambda}^{t+1} \rangle  \evv{ =}\\ &
\evv{=} - A_{t+1} \langle \nabla \varphi(\boldsymbol{\lambda}^{\evv{t+1}}), \boldsymbol{\lambda}^{t+1} \rangle  + \alpha_{t+1} \langle \nabla \varphi(\boldsymbol{\lambda}^{t+1}), \,
\boldsymbol{\lambda} \rangle + A_t \langle \nabla \varphi(\boldsymbol{\lambda}^{t+1}), \, \mathbf{y}^t \rangle \evv{ =}\\ 
&  = 
\alpha_{\evv{t}+1} 
\evv{\langle}
\nabla \varphi(\boldsymbol{\lambda}^{\evv{t+1}}), \, \boldsymbol{\lambda} - \mathbf{z}^{\evv{t}}
\evv{\rangle}.
\end{eqnarray*}
Тогда имеем
$$
A_t \left ( \varphi(\boldsymbol{\lambda}^{t+1}) + \langle \nabla \varphi(\boldsymbol{\lambda}^{t+1}),
\mathbf{y}^t - \boldsymbol{\lambda}^{t+1} \rangle \right ) +
\frac{L}{2} \| \boldsymbol{\lambda} - \mathbf{z}^t \|_2^2 
$$
$$
+
\alpha_{t+1}
\left [ 
\varphi(\boldsymbol{\lambda}^{t+1}) + 
		\evv{\langle}	\nabla \varphi (\boldsymbol{\lambda}^{t+1}), \, \boldsymbol{\lambda} - \boldsymbol{\lambda}^{t+1}
		\evv{\rangle}
			\right ] =
$$
\begin{equation}
\label{eq_last_eq_lem1}
= A_{t+1} \varphi(\boldsymbol{\lambda}^{t+1}) + 
\frac{L}{2} \| \boldsymbol{\lambda} - \mathbf{z}^t \|_2^2 +
\alpha_{t+1} 
\evv{\langle}
\nabla \varphi(\boldsymbol{\lambda}^{t+1}), \,
\boldsymbol{\lambda} - \mathbf{z}^t
\evv{\rangle}.
\end{equation}
После замены последнего выражения в~\eqref{eq_last_in_lem1} на \eqref{eq_last_eq_lem1}
можно воспользоваться расширенным вариантом
неравенства Фенхеля для сопряженных функций~\cite{Nest_dis}:
$$
\langle \mathbf{g}, \mathbf{s} \rangle + \frac{\xi}{2} \|\mathbf{s}\|^2 \, \ge \, - \frac{1}{2 \xi} \| \mathbf{g} \|_*^2, \, \mathbf{g} \, \in \, \mathbb{E}^*, \, \,  \mathbf{s} \, \in \, \mathbb{E},
$$
где
$\mathbb{E}$~--- конечномерное вещественное векторное пространство, $\mathbb{E}^*$~--- пространство линейных функций на~$\mathbb{E}$ (двойственное пространство),
норма в двойственном пространстве $\|\mathbf{g}\|_* = 
\max\limits_\mathbf{x} \{ \langle \mathbf{g}, \mathbf{x}\rangle  \, | \, \|\mathbf{x} \|_E = 1 \}$.
В нашем случае $\mathbf{g} = \nabla \varphi (\boldsymbol{\lambda}^{t+1})$,
$\mathbf{s} = \boldsymbol{\lambda} - \mathbf{z}^t$, $\xi = \frac{L}{\alpha_{t+1}}$.
Следовательно, 
\begin{equation}
\label{eq_lem1_after_fen}
\psi_{t+1}(\mathbf{z}^{t+1}) \, \ge \, 
A_{t+1} \varphi(\boldsymbol{\lambda}^{t+1}) - \frac{\alpha_{t+1}^2}{2L}
\|\nabla \varphi(\boldsymbol{\lambda}^{t+1})\|^2_2.
\end{equation}
Для завершения доказательства леммы требуется
показать, что $A_{t+1} \varphi(\mathbf{y}^{t+1})$ 
меньше, чем правая часть неравенства
в~\eqref{eq_lem1_after_fen}.

В силу $L$-гладкости функции $\varphi(\boldsymbol{\lambda})$
(см. утверждение~\ref{2_cor_sopr})
\begin{align*}
    \varphi(\mathbf{y}^{t+1}) & \, \le \,
\varphi(\boldsymbol{\lambda}^{t+1}) + \langle \nabla \varphi(\boldsymbol{\lambda}^{t+1}), \mathbf{y}^{t+1} - \boldsymbol{\lambda}^{t+1}\rangle + \frac{L}{2}
\| \mathbf{y}^{t+1} - \boldsymbol{\lambda}^{t+1} \|_2^2 = \\ &
=\min_{\boldsymbol{\lambda}}
\left \{ \varphi(\boldsymbol{\lambda}^{t+1}) + \langle \nabla \varphi (\boldsymbol{\lambda}^{t+1}) , \boldsymbol{\lambda} - \boldsymbol{\lambda}^{t+1}\rangle  + \frac{L}{2} \| \boldsymbol{\lambda} - \boldsymbol{\lambda}^{t+1} \|_2^2 
			\right
			\}\\ &  =
\varphi(\boldsymbol{\lambda}^{t+1}) - \frac{1}{2 L} \| \nabla \varphi (\boldsymbol{\lambda}^{t+1}) \|_2^2.		
\end{align*}
После умножения обеих частей полученного неравенства на $A_{t+1}$:
$$
A_{t+1} \varphi(\mathbf{y}^{t+1}) \, \le \, A_{t+1} \varphi(\boldsymbol{\lambda}^{t+1}) - \frac{A_{t+1}}{2 L} \| \nabla \varphi (\boldsymbol{\lambda}^{t+1}) \|_2^2.
$$
В силу того, что для коэффициентов БГМ $\alpha_{t+1}^2 \, \le \, A_{t+1}$, получаем
\begin{equation}
\label{eq_aphi_alam}
A_{t+1} \varphi(\mathbf{y}^{t+1}) \, \le \, 
A_{t+1} \varphi(\boldsymbol{\lambda}^{t+1}) - \frac{\alpha_{t+1}^2}{2L}
\|\nabla \varphi(\boldsymbol{\lambda}^{t+1})\|^2_2.
\end{equation}
Следовательно, в силу~\eqref{eq_lem1_after_fen} и \eqref{eq_aphi_alam}
$A_{t+1} \varphi(\mathbf{y}^{t+1}) \, \le \, \psi_{t+1}(\mathbf{z}^{t+1})$,
что и требовалось доказать.~ 

\textbf{Доказательство леммы~\ref{lem_est_1}.}
Определим следущее множество: 
$$
    \Lambda_{2\hat R} = \{\boldsymbol{\lambda} \in \mathbb{R}_{+}^m: \|\boldsymbol{\lambda}\|_2 \leq 2\hat R \}.
$$
\evv{$\hat R$ определяется в силу следующих неравенств:}
$$
||\boldsymbol{\lambda}^0-\boldsymbol{\lambda}^*||_2 + ||\boldsymbol{\lambda}^0||_2 \leq ||\boldsymbol{\lambda}^*||_2 + 2||\boldsymbol{\lambda}^0||_2 \leq 3 R = \hat R.
$$ 
При этом все $\boldsymbol{\lambda}^t$ будут принадлежать $\Lambda_{2\hat R}$, так как 
\begin{align*}
    ||\boldsymbol{\lambda}^t||_2 &\leq ||\boldsymbol{\lambda}^t-\boldsymbol{\lambda}^*||_2 + ||\boldsymbol{\lambda}^* - \boldsymbol{\lambda}^0||_2+||\boldsymbol{\lambda}^0||_2 \\ &\leq  2||\boldsymbol{\lambda}^* - \boldsymbol{\lambda}^0||_2+||\boldsymbol{\lambda}^0||_2 \evv{\, \leq \, } 2||\boldsymbol{\lambda}^*||_2 + 3||\boldsymbol{\lambda}^0||_2 \evv{\, \leq \, } 5 R \leq 2\hat R,
\end{align*}
где для второго неравенства учитывалось, что $||\boldsymbol{\lambda}^t-\boldsymbol{\lambda}^*||_2 \leq ||\boldsymbol{\lambda}^* - \boldsymbol{\lambda}^0||_2$,
\evv{$t = 0, \, 1, \, \ldots$}

\evv{Последнее неравенство можно доказать
следующим образом. 
Для любого $\boldsymbol{\lambda} \, \in \, \RR_{+}^m$ в силу
леммы~\ref{lem_fgm_pd}
и сильной выпуклости
функции $\psi_t(\boldsymbol{\lambda})$
с константой $L$ 
верно
$$
A_t \varphi(\mathbf{y}^t) + \frac{L}{2} \|  
\boldsymbol{\lambda} - \mathbf{z}^t\|_2^2
\, \le \, 
\psi_t(\mathbf{z}^t) +
 \frac{L}{2} \|  
\boldsymbol{\lambda} - \mathbf{z}^t \|_2^2
\, \le \,
\psi_t(\boldsymbol{\lambda}) 
$$
\begin{equation}
\label{ap_eq_lam_bound1}
= 
\sum\limits_{j = 0}^t \alpha_j \left [ 
\varphi(\boldsymbol{\lambda}^j) + \langle
			\nabla \varphi (\boldsymbol{\lambda}^j), \boldsymbol{\lambda} - \boldsymbol{\lambda}^j\rangle 
			\right ] 
			+
 \frac{L}{2} \|  
\boldsymbol{\lambda} - \boldsymbol{\lambda}^0 \|_2^2.
\end{equation}
Последнее выражение в~\eqref{ap_eq_lam_bound1} 
в силу выпуклости функции $\varphi(\boldsymbol{\lambda})$
можно
оценить сверху как 
$A_t \varphi (\boldsymbol{\lambda}) +
 \frac{L}{2} \|  
\boldsymbol{\lambda} - \boldsymbol{\lambda}^0 \|_2^2$.
Тогда при $\boldsymbol{\lambda} =
\boldsymbol{\lambda}^*$
$$
\frac{L}{2} \|  
\boldsymbol{\lambda}^* -
\mathbf{z}^t\|_2^2 
\, \le \, 
A_t \left ( 
\varphi(\mathbf{y}^t) - \varphi (\boldsymbol{\lambda}^*) 
\right )
+ \frac{L}{2} \|  
\boldsymbol{\lambda}^* -
\mathbf{z}^t\|_2^2
\, \le \, 
\frac{L}{2} \|  
\boldsymbol{\lambda}^* - \boldsymbol{\lambda}^0 \|_2^2.
$$
Следовательно,
\begin{equation}
\label{ap_eq_lam_bound2}
    \|  
\boldsymbol{\lambda}^* -
\mathbf{z}^t\|_2
\, \le \,
\|  
\boldsymbol{\lambda}^* - \boldsymbol{\lambda}^0 \|_2.
\end{equation}
Поскольку $\mathbf{y}^t$
в алгоритме~\ref{alg_fgm}
определяется с помощью
шага метода проекции градиента
для выпуклой функции
$\varphi(\boldsymbol{\lambda})$,
последовательность генерируемых
алгоритмом точек
$\mathbf{y}^t$, $t = 0, \, 1, \, \ldots$ также будет
ограничена (доказательство этого
факта см., например, в~\cite[лемма~9.17, стр.~183]{beck_matlab}
или в~\cite[стр.~265]{bubeck}):
\begin{equation}
\label{ap_eq_lam_bound3}
    \|  
\boldsymbol{\lambda}^* -
\mathbf{y}^t\|_2
\, \le \,
\|  
\boldsymbol{\lambda}^* - \boldsymbol{\lambda}^0 \|_2.
\end{equation}
Далее
$$
\|  
\boldsymbol{\lambda}^{t+1} -
\boldsymbol{\lambda}^{*}\|_2
=
\|  
\tau_t ( \mathbf{z}^t -
\boldsymbol{\lambda}^{*}) +
(1 - \tau_t) 
(\mathbf{y}^t -
\boldsymbol{\lambda}^{*})
\|_2
\, \le \,
\tau_t \| \mathbf{z}^t -
\boldsymbol{\lambda}^{*} \|_2
+
(1 - \tau_t) 
\| \mathbf{y}^t -
\boldsymbol{\lambda}^{*})
\|_2.
$$
Из последнего неравенства
с помощью~\eqref{ap_eq_lam_bound2}, \eqref{ap_eq_lam_bound3}
получаем нужный результат:
$$
\|  
\boldsymbol{\lambda}^{t+1} -
\boldsymbol{\lambda}^{*}\|_2
\, \le, 
\|  
\boldsymbol{\lambda}^{*} -
\boldsymbol{\lambda}^{0}\|_2,
\, t = -1, \, 0, \, 1, \, \ldots
$$
}.


В силу леммы~\ref{lem_fgm_pd}
\begin{eqnarray*}
A_N \varphi(\mathbf{y}^N) &  \le & \min_{\boldsymbol{\lambda} \, \in \, 
\RR_+^m} 
\left \{
\frac{L}{2} \|\boldsymbol{\lambda} - \boldsymbol{\lambda}^0\|_2^2 +
\sum\limits_{t = 0}^N \alpha_t \left [ 
\varphi(\boldsymbol{\lambda}^t) + 
			\langle 
			\nabla \varphi (\boldsymbol{\lambda}^t), \boldsymbol{\lambda} - \boldsymbol{\lambda}^t \rangle
			\right ] 
\right \} \, \\
& \le &
\min_{\boldsymbol{\lambda} \, \in \, 
\Lambda_{2\hat R}} 
\left \{
\frac{L}{2} \|\boldsymbol{\lambda} - \boldsymbol{\lambda}^0\|_2^2 +
\sum\limits_{t = 0}^N \alpha_t \left [ 
\varphi(\boldsymbol{\lambda}^t) + 
			\langle 
			\nabla \varphi (\boldsymbol{\lambda}^t), \boldsymbol{\lambda} - \boldsymbol{\lambda}^t \rangle
			\right ] 
\right \} \,\\ &
\overset{\circledOne}{\le} &
 \min_{\boldsymbol{\lambda} \, \in \, 
\Lambda_{2\hat R}}
\left \{
\sum\limits_{t = 0}^N \alpha_t \left [ 
\varphi(\boldsymbol{\lambda}^t) + \langle 
			\nabla \varphi (\boldsymbol{\lambda}^t), \boldsymbol{\lambda} - \boldsymbol{\lambda}^t \rangle
			\right ] 
\right \} +  \frac{ 37L\hat R^2}{9},
\end{eqnarray*}
где $\circledOne$ выполняется, так как
\begin{equation}
\label{eq_fgm_ll0}
\|\boldsymbol{\lambda} - \boldsymbol{\lambda}^0\|^2 \, \le \,
2 \|\boldsymbol{\lambda}\|^2 + 2 \|\boldsymbol{\lambda}^0\|^2 \, \le \, 
8 \hat R^2 + \frac{2}{9} \hat R^2 = \frac{74}{9} \hat R^2. 
\end{equation}
После подстановки определений двойственной целевой
функции $\varphi(\boldsymbol{\lambda}^t)$~\eqref{eq_phi_first} и ее градиента~$\nabla \varphi(\boldsymbol{\lambda}^t)$ (см. утверждение~\ref{2_cor_dem}):
\begin{align*}
\sum\limits_{t = 0}^N & \alpha_t \left [ 
\varphi(\boldsymbol{\lambda}^t) + 
			\langle 
			\nabla \varphi (\boldsymbol{\lambda}^t), \boldsymbol{\lambda} - \boldsymbol{\lambda}^t \rangle
			\right ] \\ & = 
\sum\limits_{t = 0}^N \alpha_t
\Bigg (
\langle \boldsymbol{\lambda}^t, \mathbf{b} \rangle + \sum_{k = 1}^n  \left ( u_k(x^t_k(\boldsymbol{\lambda}^t)) - \langle \boldsymbol{\lambda}^t, \mathbf{C}_k x_k^t(\boldsymbol{\lambda}^t) \rangle\right )
\\
& + \langle\mathbf{b} - \sum_{k = 1}^n \mathbf{C}_k x_k^t(\boldsymbol{\lambda}^t) ,
 \boldsymbol{\lambda} - \boldsymbol{\lambda}^t \rangle
\Bigg )  \\ &=   
\sum\limits_{t = 0}^N \alpha_t
\left (
\sum_{k = 1}^n u_k(x_k^t(\boldsymbol{\lambda}^t)) +
 \langle \boldsymbol{\lambda},  \mathbf{b} - \sum_{\evv{k} = 1}^n \mathbf{C}_\evv{k} x_k^t(\boldsymbol{\lambda}^t) \rangle
\right )
 \\ &  \le \,
A_N \left (
U(\mathbf{\hat{x}}^{N}) + \langle \boldsymbol{\lambda} , \mathbf{b} -  C\mathbf{\hat{x}}^N \rangle
\right ), 
\end{align*}
где последнее неравенство выполняется в силу вогнутости функций полезности.

Итак, 
\begin{align*}
A_N \varphi(\mathbf{y}^N) \, & \le \, A_N U(\mathbf{\hat{x}}^{N}) +
\frac{ 37L\hat R^2}{9} +
A_N \min_{\boldsymbol{\lambda} \, \in \, 
\Lambda_{2 \hat R}} 
\left \{ \langle
\boldsymbol{\lambda}, \mathbf{b} -  C\mathbf{\hat{x}}^N \rangle
\right \} \\ &
= A_N U(\boldsymbol{\hat{x}}^{N}) +
\frac{ 37L\hat R^2}{9} -
A_N \max_{\boldsymbol{\lambda} \, \in \, 
\Lambda_{2 \hat R}} 
\left \{\langle
\boldsymbol{\lambda}, C \mathbf{\hat{x}}^N - \mathbf{b} \rangle
\right \} \\
& = A_N U(\boldsymbol{\hat{x}}^{N}) +
\frac{ 37L\hat R^2}{9}- 2 \hat R \pd{A_N} \left \|\left ( C \mathbf{\hat{x}}^N - \mathbf{b}\right )_+ \right \|_2.
\end{align*}
Из этого получаем оценку \eqref{eq_aur37}. ~

%
\subsection*{Доказательство леммы~\ref{lem_stoch_main}}

Для доказательства леммы~\ref{lem_stoch_main} сначала приведем доказательства нескольких вспомогательных технических лемм.
\begin{Lem}\label{lem:new_recurrence_lemma_appendix}
     Пусть $A, B$ и $\{r_t\}_{t=0}^N$~--- неотрицательные числа, такие, что для любого $l = 1, \, \ldots, \, N$ выполняется неравенство
    \begin{equation}\label{eq:new_bound_for_r_l_appendix}
         \frac{1}{2} r_l^2 \, \le \, A r_0^2 + B r_0\sqrt{\sum\limits_{t=0}^{l-1} r_t^2}.
     \end{equation}
     Тогда верно следующее неравенство:
     \begin{equation}\label{eq:new_recurrence_lemma_appendix}
         r_l \, \le \, C r_0,
     \end{equation}
     где $C$~--- положительная константа, для которой выполняется $C^2 \, \ge \, \max \left \{1, \, 2A + 2BC \sqrt{N} \right \}$, т.е., в частности, можно выбрать 
     $$
     C = \max \left \{1, \, B\sqrt{N} + \sqrt{B^2N + 2A} \right \}.
     $$
\end{Lem}
\textbf{Доказательство.}
  Докажем \eqref{eq:new_recurrence_lemma_appendix} по индукции. Для $l=0$ неравенство выполнено, так как $C \, \ge \, 1$. Пусть \eqref{eq:new_recurrence_lemma_appendix} выполнено для всех $l < N$. Докажем,
  что оно выполнено и для $l+1$. Действительно,
     \begin{eqnarray*}
        r_{l+1} \,  \overset{\eqref{eq:new_bound_for_r_l_appendix}}{\le} \, 
        \sqrt{2}\sqrt{Ar_0^2 + B r_0\sqrt{\sum\limits_{t=0}^{l}r_t^2}}
    \,    \overset{\eqref{eq:new_recurrence_lemma_appendix}}{\le} \,
        r_0\sqrt{2}\sqrt{A + B C \sqrt{N}}
    \\
        = 
        r_0\underbrace{\sqrt{2A + 2 B C \sqrt{N}}}_{\le \, C} \, \le \, Cr_0. 
     \end{eqnarray*}

\begin{Lem}
\label{lem_tail_est}
Пусть для последовательностей неотрицательных коэффициентов $\{R_t\}_{t \, \ge \, 0}$ и случайных векторов $\{ \boldsymbol{\eta}^t\}_{t \, \ge \, 0}$,  $\{\mathbf{a}^t\}_{t \, \ge \, 0}$  для всех $l = 1, \, \ldots, \, N$ выполняется неравенство
   \begin{eqnarray}
        \frac{1}{2}R_l^2 \le A + u\sum\limits_{t=0}^{l-1}\langle \boldsymbol{\eta}^t, \, \mathbf{a}^t\rangle,
        \label{eq:radius_recurrence_appendix}
    \end{eqnarray}
    где  $A$~--- неотрицательная константа, $d \, \ge \, 1$~--- положительная константа,
    $\|\mathbf{a}^t\|_2 \, \le \, \widetilde{R}_t d$ и $\widetilde{R}_t = \max \left \{\widetilde{R}_{t-1}, \, R_t \right \}$ для всех $t \, \ge \, 1$, $\widetilde{R}_0 = R_0$, $\widetilde{R}_t$ зависит только от  $\boldsymbol{\eta}^0, \, \ldots, \, \boldsymbol{\eta}^t$. Пусть также вектор $\mathbf{a}^t$~\ev{---} это функция от  $\boldsymbol{\eta}^0, \, \ldots, \, \boldsymbol{\eta}^{t-1}$  $\forall \, t \, \ge \, 1$, $a^0$~\ev{---} постоянный вектор  и для любого $t \, \ge \, 0$ 
    \begin{eqnarray*}
\EE\left[\boldsymbol{\eta}^t\mid \{\boldsymbol{\eta}^k\}_{k=0}^{t-1}\right] = 0,\quad 
\EE\left[\exp\left({\|\boldsymbol{\eta}^t\|_2^2}{\sigma^{-2}}\right)\mid\{\boldsymbol{\eta}^k\}_{k=0}^{t-1}\right] \le \exp(1). \label{eq:eta_k_properties_appendix}
    \end{eqnarray*}
\ev{Т}огда с вероятностью $1 - 2\delta$ выполняются  следующие неравенства:
  \begin{eqnarray*}
        \widetilde{R}_l \le JR_0 \quad \text{ и }\label{eq:tails_estimate_radius_appendix} \quad
        A + u\sum\limits_{t=0}^{l-1}\langle \boldsymbol{\eta}^t, \mathbf{a}^t \rangle \le A + udD\sqrt{\sigma^2 g(N)NJ}\widetilde{R}_0^2
  \end{eqnarray*}
$\forall \, \, l=1, \, \ldots, \, N$ одновременно, где $D$~\ev{---} положительная константа, 
$$
\ev{F} = 2\sigma^2d^2 N (2ud)^{N}\left(2A + ud \widetilde{R}_0^2 + 12 ud \ln\frac{N}{\delta} \sigma^2 N \right),
$$
$\ev{f}=d^2\sigma^2\widetilde{R}_0^2$, $g(N) = \ln\left(\frac{N}{\delta}\right) + \ln\ln\left(\frac{\ev{F}}{\ev{f}}\right)$ и 
$$
J = \max\left\{1, \, \frac{1}{\widetilde{R}_0}udD\sqrt{\sigma^2 g(N)} + \sqrt{\frac{1}{\widetilde{R}_0^2}u^2d^2 C_1^2\sigma^2g(N) + \frac{2A}{R_0^2} }\right\}.$$
\end{Lem}
\textbf{Доказательство.}
Применим ко второму слагаемому из правой части \eqref{eq:radius_recurrence_appendix} неравенство Коши-Буняковского\ev{:}
  \begin{eqnarray}\label{eq:radius_recurrence2}
        \frac{1}{2}R_l^2 \, \le \,  A + ud\sum\limits_{t=0}^{l-1}\|\boldsymbol{\eta}^t\|_2\widetilde{R}_t \, \le \,  A + \frac{ud}{2}\sum\limits_{t=0}^{l-1}\widetilde{R}_t^2 + \frac{ud}{2}\sum\limits_{t=0}^{l-1}\|\boldsymbol{\eta}^t\|_2^2.
    \end{eqnarray}
    
По теореме~2.1 из~\cite{JudNem}
\evv{
\begin{multline}
\label{app_eq_JudNem}
(\forall \, N \, \ge \, 1, \, \forall
\, \gamma \, \ge \, 0): \quad
\mathbb{P} 
\left \{ 
\left \| \sum\limits_{t=0}^{N-1} \boldsymbol{\eta}^t \right \|_2 \, \ge \,
(\sqrt{2} + \sqrt{2} \gamma) \sqrt{ \sum\limits_{t=0}^{N-1} \sigma_t^2} 
\right \} \\
\le \,
\exp \left(-\frac{\gamma^2}{3} \right ).
\end{multline}
Тогда с вероятностью не меньшей, чем
\begin{equation}
\label{app_eq_gamma}
1 - \frac{\delta}{N} = 1 - \exp \left(-\frac{\gamma^2}{3} \right )
\end{equation}
}
выполняется следующее неравенство 
\begin{eqnarray}
\label{app_eq_eta_t}
        \|\boldsymbol{\eta}^t\|_2 \le \sqrt{2}\left(1 + \sqrt{3\ln\frac{N}{\delta}}\right)\sigma \le  2\sqrt{6\ln\frac{N}{\delta}}\sigma.
\end{eqnarray}
\evv{Действительно, выражая из~\eqref{app_eq_gamma} $\gamma$,
получаем, что $\gamma = \sqrt{3\ln\frac{N}{\delta}}$.
Поcле подстановки данного
выражения в~\eqref{app_eq_JudNem}
и выбора единого $\sigma \, \in \, \RR_+$
вместо последовательности $\sigma_t$, $t = 0, \, \ldots, \, N-1$, получаем оценку~\eqref{app_eq_eta_t}.}

Объединяя полученные неравенства, получаем\ev{,} что с вероятностью\ev{, большей или равной}~$1 - \delta$ неравенство 
\begin{eqnarray*}
        \frac{1}{2}R_l^2  \le  A + \frac{ud}{2}\sum\limits_{t=0}^{l-1}\widetilde{R}_t^2 + 12 ud \ln\frac{N}{\delta} \sigma^2 l
    \end{eqnarray*}
выполняется для всех $l = 1, \, \ldots, \, N$ одновременно. Заметим, что последнее слагаемое в полученной оценк\ev{е~---}  неубывающая функция от $l$. Определим $\hat l$ как наибольшее целое число\ev{,} для которого выполнено $\hat l \leq l$ и $\widetilde{R}_{\hat{l}} = R_{\hat{l}}$. Тогда получаем, что $R_{\hat{l}} = \widetilde{R}_{\hat{l}} = \widetilde{R}_{\hat{l}+1} = \ldots = \widetilde{R}_{l}$ и\ev{,} следовательно\ev{,} с вероятностью $\ge 1 - \delta$
    \begin{eqnarray*}
        \frac{1}{2}\widetilde{R}_l^2 \le A + \frac{ud}{2}\sum\limits_{t=0}^{\hat{l}-1}\widetilde{R}_t^2 + 12 ud \ln\frac{N}{\delta} \sigma^2 \hat{l}
        \le  A + \frac{ud}{2}\sum\limits_{t=0}^{l-1}\widetilde{R}_t^2 + 12 ud \ln\frac{N}{\delta} \sigma^2 l 
        \\
        \forall \, \, l=1, \, \ldots, \, N.
    \end{eqnarray*}
Получаем\ev{,} что с вероятностью $\ge 1 - \delta$ верна следующая оценка: 
 \begin{eqnarray*}
        \widetilde{R}_l^2 &\le& 2A + ud\sum\limits_{t=0}^{l-1}\widetilde{R}_k^2 + 24 ud \ln\frac{N}{\delta} \sigma^2 l \\ 
        &\le& 2A\underbrace{(1+ud)}_{\le 2ud} + \underbrace{(ud + u^2d^2)}_{\le 2u^2d^2}\sum\limits_{t=0}^{l-2}\widetilde{R}_t^2 + 24 ud \ln\frac{N}{\delta} \sigma^2\underbrace{( l+ud(l-1))}_{\leq 2udl} \notag\\
        &\le& 2ud\left(2A + ud\sum\limits_{t=0}^{l-2}\widetilde{R}_t^2 + 24 ud \ln\frac{N}{\delta} \sigma^2 l \right) \quad \forall \, \, l = 1, \, \ldots, \, N.\notag
    \end{eqnarray*}
Применяя данную оценку рекурс\ev{и}вно\ev{,} получаем, что с вероятностью $\ge 1 - \delta$ верно
$$
\widetilde{R}_l^2 \le (2ud)^{l}\left(2A + ud \widetilde{R}_0^2 + 24 ud \ln\frac{N}{\delta} \sigma^2 l \right).
$$
Далее рассмотрим последовательность случайных вел\ev{и}чин $\xi^t = \langle \boldsymbol{\eta}^t, \, \mathbf{a}^t \rangle$. Зам\ev{e}тим, что  $\EE\left[\xi^t\mid \xi^0, \, \ldots, \, \xi^{t-1}\right] = \left\langle \EE\left[\boldsymbol{\eta}^t\mid \boldsymbol{\eta}^0, \, \ldots, \, \boldsymbol{\eta}^{k-1}\right], \, \mathbf{a}^t \right\rangle = 0$\ev{,} тогда\ev{,} используя неравенство Коши-Буняковского\ev{,} получаем, что 
 \begin{eqnarray*}
        \EE\left[\exp\left(\frac{(\xi^t)^2}{\sigma^2d^2\widetilde{R}_t^2}\right)\mid \xi^0, \, \ldots, \, \xi^{t-1}\right] &\le& \EE\left[\exp\left(\frac{\|\boldsymbol{\eta}^t\|_2^2 d^2\widetilde{R}_t^2}{\sigma^2 d^2\widetilde{R}_t^2}\right)\mid \boldsymbol{\eta}^0, \, \ldots, \, \boldsymbol{\eta}^{t-1}\right]\\
        &=& \EE\left[\exp\left(\frac{\|\boldsymbol{\eta}^t\|_2^2}{\sigma^2}\right)\mid \boldsymbol{\eta}^0, \, \ldots, \, \boldsymbol{\eta}^{t-1}\right] \\
        \le \, \exp(1).
    \end{eqnarray*}
Определим $\hat\sigma_t^2 = \sigma^2d^2\widetilde{R}_t^2$, тогда с вероятностью $\ge 1 - \delta$ выполняется 
\begin{eqnarray*}
 \sum\limits_{t=0}^{l-1} \hat\sigma_t^2 & \leq  &\sigma^2d^2 l (2ud)^{l}\left(2A + ud \widetilde{R}_0^2 + 24 ud \ln\frac{N}{\delta} \sigma^2 l \right) \\ & \leq& \sigma^2d^2 N (2ud)^{N}\left(2A + ud \widetilde{R}_0^2 + 24 ud \ln\frac{N}{\delta} \sigma^2 N \right) := \dfrac{\ev{F}}{2}
 \end{eqnarray*}
для всех  $l = 1, \, \ldots, \, N$ одновременно, где 
$$
\ev{F} = 2\sigma^2d^2 N (2ud)^{N}\left(2A + ud \widetilde{R}_0^2 + 24 ud \ln\frac{N}{\delta} \sigma^2 N \right).
$$ 

Используя следствие~8 из \cite{jin2019short} для $b=\hat\sigma_0^2$\ev{,} получаем для любого $l=1, \, \ldots, \, N$ с вероятностью $\ge 1-\frac{\delta}{N}$ следующую оценку:
    \begin{equation}
    \label{ap_eq_libo}
        \text{либо } \sum\limits_{t=0}^{l-1}\hat\sigma_t^2 \ge \ev{F}, \text{ либо } \left|\sum\limits_{t=0}^{l-1}\xi^t\right| \le C_1\sqrt{\sum\limits_{t=0}^{l-1}\hat\sigma_t^2\left(\ln\left(\frac{N}{\delta}\right) + \ln\ln\left(\frac{\ev{F}}{\ev{f}}\right)\right)},
    \end{equation}
где $C_1 >0$~\ev{---} константа, которая не зависит от $\ev{F}$ и $\ev{f}$. 

Далее\ev{,} объединяя полученные оценки\ev{,} получаем, что
\evv{оценка~\eqref{ap_eq_libo}}
с вероятностью $\ge 1-\delta$ 
\evv{верна } 
для всех  $l=1, \, \ldots, \, N$ одновременно. 
 
 Учитывая выбор $\ev{F}$\ev{,} получаем, что с вероятностью $\ge 1-2\delta$
\begin{equation*}
\left|\sum\limits_{t=0}^{l-1}\xi^t\right| \le C_1\sqrt{\sum\limits_{t=0}^{l-1}\hat\sigma_t^2\left(\ln\left(\frac{N}{\delta}\right) + \ln\ln\left(\frac{\ev{F}}{\ev{f}}\right)\right)}
    \end{equation*}
    для всех  $l=1, \, \ldots, \, N$ одновременно. 
    
    Для удобства дальнейших рассуждений обозначим $g(N) := \ln\left(\frac{N}{\delta}\right) + \ln\ln\left(\frac{\ev{F}}{\ev{f}}\right) \approx \ln\left(\frac{N}{\delta}\right) $, пр\ev{е}небрегая константой. Использу\ev{я} $\hat\sigma_t^2 = \sigma^2d^2\widetilde{R}_t^2$\ev{,}  получаем, что с вероятностью $\ge 1-2\delta$ справедлива следующая оценка\ev{:} 
 \begin{eqnarray}
        \frac{1}{2}\widetilde{R}_l^2 &\le& A + u\sum\limits_{t=0}^{l-1}\underbrace{\langle \boldsymbol{\eta}^t, \, \mathbf{a}^t\rangle }_{\xi^t}
    \,    \le \, A + udD\sqrt{\sigma^2 g(N)}\sqrt{\sum\limits_{t=0}^{l-1}\widetilde{R}_t^2}\notag\\
\end{eqnarray}
для всех  $l=1, \, \ldots, \, N$ одновременно. Выбирая в качестве $A=\frac{A}{\widetilde{R}_0^2}$\ev{,}  $B =\frac{1}{\widetilde{R}_0}udC_1\sqrt{\sigma^2 g(N)}$, $r_t = \widetilde{R}_t $\ev{,} из леммы~\ref{lem:new_recurrence_lemma_appendix} получаем\ev{,} что с вероятностью  $1-2\delta$ выполняется 
    \begin{eqnarray*}
        \widetilde{R}_l \le JR_0
    \end{eqnarray*}
    для всех $l=1, \, \ldots, \, N$ одновременно, где  
    $$
    J = \max\left\{1, \, \frac{1}{\widetilde{R}_0}udC_1\sqrt{\sigma^2 g(N)} + \sqrt{\frac{1}{\widetilde{R}_0^2}u^2d^2C_1^2\sigma^2g(N) + \frac{2A}{R_0^2} }\right\}.
    $$ 
    Отсюда получаем, что с вероятностью $1-2\delta$ оценка
 \begin{eqnarray*}
         A + u\sum\limits_{t=0}^{l-1}\langle \boldsymbol{\eta}^t, \mathbf{a}^t\rangle
        \le A + udC_1\sqrt{\sigma^2 g(N)lJ}\widetilde{R}_0^2 \le A + udC_1\sqrt{\sigma^2 g(N)NJ}\widetilde{R}_0^2\notag\\
\end{eqnarray*}
\ev{верна} для всех  $l=1, \, \ldots, \, N$ одновременно. 
~

\textbf{Доказательство леммы~\ref{lem_stoch_main}}

Для $\boldsymbol{\lambda} \, \in \, \RR_+^m$ 
$$
\| \boldsymbol{\lambda}^{t+1} - \boldsymbol{\lambda} \|_2^2 =
\|
[\boldsymbol{\lambda}^t - \beta \nabla \varphi(\boldsymbol{\lambda}^t, \, \xi^t)]_+ - \boldsymbol{\lambda} \|_2^2 \, \le \,
\| \boldsymbol{\lambda}^t - \boldsymbol{\lambda} \|_2^2 - 
2 \beta \langle \nabla \varphi(\boldsymbol{\lambda}^t, \, \xi^t),
\boldsymbol{\lambda}^t - \boldsymbol{\lambda}\rangle 
$$
$$
+ 
\beta^2 
\| \nabla \varphi(\boldsymbol{\lambda}^t, \, \xi^t) \|_2^2,
$$
т.е. 
\begin{equation}
\label{eq_stoch_t1}
0 \, \le \, \frac{1}{2 \beta}
\left (
\| \boldsymbol{\lambda}^t - \boldsymbol{\lambda} \|_2^2 -
\| \boldsymbol{\lambda}^{t+1} - \boldsymbol{\lambda} \|_2^2 
\right ) +
\langle \nabla \varphi(\boldsymbol{\lambda}^t, \, \xi^t),
\boldsymbol{\lambda} - \boldsymbol{\lambda}^t \rangle +
\frac{\beta}{2} 
\| \nabla \varphi(\boldsymbol{\lambda}^t, \, \xi^t) \|_2^2.    
\end{equation}
После прибавления к обеим сторонам неравенства~\eqref{eq_stoch_t1} $\varphi(\boldsymbol{\lambda}^t)$, умножения на $N$
и суммирования от $0$ до $N-1$:
\begin{align}
\dfrac{1}{N}\sum_{t = 0}^{N-1} \varphi(\boldsymbol{\lambda}^t) & \, \le \,
\dfrac{1}{N} \sum_{t = 0}^{N-1} 
\Bigg \{
\varphi(\boldsymbol{\lambda}^t) +
\langle \nabla \varphi(\boldsymbol{\lambda}^t, \, \xi^t), \,
\boldsymbol{\lambda} - \boldsymbol{\lambda}^t\rangle +
\frac{\beta}{2} 
\| \nabla \varphi(\boldsymbol{\lambda}^t, \, \xi^t) \|_2^2 \notag \\ &
+ \frac{1}{2 \beta}
\Big (
\| \boldsymbol{\lambda}^t - \boldsymbol{\lambda} \|_2^2 -
\| \boldsymbol{\lambda}^{t+1} - \boldsymbol{\lambda} \|^2_2
\Big )
\Bigg \}.
\end{align}
В силу выпуклости $\varphi(\boldsymbol{\lambda})$ 
для $\boldsymbol{\hat{\lambda}}^N = \dfrac{1}{N}\sum\limits_{t=0}^{N-1}
\boldsymbol{\lambda}^t$ получаем, что 
\begin{align}
\label{eq:1}
N  \varphi(\ev{\hat{\boldsymbol{\lambda}}}^N) & \, \le \,
 \sum_{t = 0}^{N-1} 
\left \{
\varphi(\boldsymbol{\lambda}^t) +
\langle \nabla \varphi(\boldsymbol{\lambda}^t, \, \xi^t),
\boldsymbol{\lambda} - \boldsymbol{\lambda}^t\rangle \right \} \notag
\\ &
+
\frac{\beta}{2} 
\| \nabla \varphi(\boldsymbol{\lambda}^t, \, \xi^t) \|_2^2 
+ \frac{1}{2 \beta}
\Big (
\| \boldsymbol{\lambda}^0 - \boldsymbol{\lambda} \|_2^2 -
\| \boldsymbol{\lambda}^{N} - \boldsymbol{\lambda} \|^2_2
\Big ).
\end{align}
Выбираем $ \boldsymbol{\lambda} =  \boldsymbol{\lambda}^*$ и прибавляем и вычитаем справа $ \sum_{t = 0}^{N-1} \langle \nabla \varphi(\boldsymbol{\lambda}^t), \,
\boldsymbol{\lambda}^* - \boldsymbol{\lambda}^t\rangle$. \ev{П}олучаем 
\begin{eqnarray}
  N  \varphi(\ev{\hat{\boldsymbol{\lambda}}}^N) &  \le &
 \sum\limits_{t = 0}^{N-1} 
\left \{
\varphi(\boldsymbol{\lambda}^t) +
\langle \nabla \varphi(\boldsymbol{\lambda}^t),
\boldsymbol{\lambda}^* - \boldsymbol{\lambda}^t\rangle \right \} +
\frac{\beta}{2} 
\| \nabla \varphi(\boldsymbol{\lambda}^t, \, \xi^t) \|_2^2 \notag \\ 
& + &  \sum\limits_{t = 0}^{N-1} \langle \nabla \varphi( \boldsymbol{\lambda}^t, \, \xi^t) - \nabla \varphi(\boldsymbol{\lambda}^t),
\boldsymbol{\lambda}^* - \boldsymbol{\lambda}^t\rangle
\notag \\
& + & \frac{1}{2 \beta}
\Big (
\| \boldsymbol{\lambda}^0 - \boldsymbol{\lambda}^* \|_2^2 -
\| \boldsymbol{\lambda}^{N} - \boldsymbol{\lambda}^* \|^2_2
\Big ).      \label{eq:3}
\end{eqnarray}
Из выпуклости $\varphi(\boldsymbol{\lambda})$ \ev{имеем}
\begin{eqnarray*}
 \sum\limits_{t = 0}^{N-1} 
\left \{
\varphi(\boldsymbol{\lambda}^t) +
\langle \nabla \varphi(\boldsymbol{\lambda}^t), \,
\boldsymbol{\lambda}^* - \boldsymbol{\lambda}^t\rangle \right \}\, \le \, \sum\limits_{t = 0}^{N-1} 
\left \{
\varphi(\boldsymbol{\lambda}^t) + \varphi(\boldsymbol{\lambda}^*) - \varphi(\boldsymbol{\lambda}^t) \right \} 
\\
\leq
\sum\limits_{t = 0}^{N-1}  \varphi(\ev{\boldsymbol{\lambda}^*}) \leq N\varphi(\boldsymbol{\lambda}^*).
\end{eqnarray*}
Подставляя полученную оценку в \eqref{eq:3}, получаем 
\begin{eqnarray}
 \frac{1}{2 \beta} \| \boldsymbol{\lambda}^{N} - \boldsymbol{\lambda}^* \|^2_2 \leq 
  \frac{1}{2 \beta} \| \boldsymbol{\lambda}^{0} - \boldsymbol{\lambda}^* \|^2_2 +
 \sum\limits_{t = 0}^{N-1} 
\langle \nabla \varphi( \boldsymbol{\lambda}^t, \, \xi^t) - \nabla \varphi(\boldsymbol{\lambda}^t), \,
\boldsymbol{\lambda}^* - \boldsymbol{\lambda}^t\rangle 
\notag \\
+
\frac{\beta}{2} 
\| \nabla \varphi(\boldsymbol{\lambda}^t, \, \xi^t) \|_2^2.   \label{eq:2}
\end{eqnarray}
Определим $R_{t}=\| \boldsymbol{\lambda}^{t} - \boldsymbol{\lambda}^* \|_2$ и $\widetilde{R}_t = \max\{\widetilde{R}_{t-1}, \, R_t\}$, причем $R_0=\widetilde{R}_0$ и\ev{,} так как\ev{,} $ \boldsymbol{\lambda}^{0}=\mathbf{0}$  и $\| \boldsymbol{\lambda}^{*}\|_2 \leq R$, то $R_0=R$. При этом по построению получаем, что $\boldsymbol{\lambda}^{t} \, \in \, B_{\widetilde{R}_t}(\boldsymbol{\lambda}^{*})$. Так же определим $\|\mathbf{a}^t \|_2= \| \boldsymbol{\lambda}^{t} - \boldsymbol{\lambda}^* \|_2 \leq \widetilde{R}_t$. Тогда \eqref{eq:2} можно переписать в следующем виде\ev{:}
\begin{eqnarray*}
 \frac{1}{2 \beta}\widetilde{R}_N^2 \leq 
  \frac{1}{2 \beta} \widetilde{R}_0^2 +
 \sum\limits_{t = 0}^{N-1} 
\langle \nabla \varphi( \boldsymbol{\lambda}^t, \, \xi^t)- \nabla \varphi(\boldsymbol{\lambda}^t), \,
\mathbf{a}^t \rangle  +
\frac{\beta}{2} 
\| \nabla \varphi(\boldsymbol{\lambda}^t, \, \xi^t) \|_2^2.   
\end{eqnarray*}
Обозначим $\boldsymbol{\eta}^t = \nabla \varphi(\boldsymbol{\lambda}^t, \, \xi^t) -
\nabla \varphi(\boldsymbol{\lambda}^t)$.
По теореме~2.1 из~\cite{JudNem}
\begin{eqnarray}
\label{estimation_1}
\mathbb{P} 
\left \{ \left
\| \sum\limits_{t=0}^{N-1} \boldsymbol{\eta}^t \right \|_2 \, \ge \,
(\sqrt{2} + \sqrt{2} \gamma) \sqrt{ \sum\limits_{t=0}^{N-1} \sigma_t^2} \, | \,
\{\xi^t\}_{t=0}^{N-1} 
\right \} \, \le \,
\exp \left(-\frac{\gamma^2}{3} \right ).
\end{eqnarray}
Используя лемму 2 из \cite{jin2019short}, получаем\ev{, что}  $$
\EE\left[\exp\left({\frac{\|\boldsymbol{\eta}^t\|_2^2}{\sigma^2}}\right)|\{\xi^k\}_{k=0}^{t-1}\right] \le \exp(1),
$$
при этом 
 $\ev{\mathbf{\eta}}^t$ зависит только от $\xi^{t-1}, \, \ldots, \, \xi^0$. Используя новые обозначения и \eqref{eq_stoch_gradM}, \ev{имеем}
\begin{eqnarray*}
\widetilde{R}_N^2 \leq 
 \widetilde{R}_0^2 + 2 \beta
 \sum\limits_{t = 0}^{N-1} 
\langle\boldsymbol{\eta}^t, \,
\mathbf{a}^t \rangle  +
\beta^2 M^2.   
\end{eqnarray*}
Тогда из леммы \ref{lem_tail_est} с константами $A= \widetilde{R}_0^2 + \beta^2 M^2$, $d=1$ и $u =  \beta$, получаем, что с вероя\ev{т}ностью $1-2\delta$,
 \ev{ где~$\frac{\delta}{N} = \exp \left(-\frac{\gamma^2}{3} \right )$,} верна следующая оценка\ev{:} 
 \begin{eqnarray}
  \label{est_part_2}
        \widetilde{R}_l \le JR_0 \quad \text{ и }
        \quad
         \sum\limits_{t=0}^{l-1}\langle \boldsymbol{\eta}^t, \, \mathbf{a}^t \rangle \, \le     \, D\sqrt{\sigma^2 g(N)NJ}\widetilde{R}_0^2
  \end{eqnarray}
$\forall \, l=1, \, \ldots, \, N$ одновременно, где $D$~\ev{---} положительная константа,
$$
\ev{F} = 2\sigma^2 N (2 \beta)^{N}\left(2A +  \beta \widetilde{R}_0^2 + 24 \ln\frac{N}{\delta} \beta \sigma^2 N \right),
$$
$\ev{f}=\sigma^2\widetilde{R}_0^2$, $g(N) = \ln\left(\frac{N}{\delta}\right) + \ln\ln\left(\frac{\ev{F}}{\ev{f}}\right)$ и 
$$
J = \max\left\{1, \, \frac{1}{\widetilde{R}_0} \beta C_1\sqrt{\sigma^2 g(N)} + \sqrt{\frac{1}{\widetilde{R}_0^2} \beta^2 C_1^2\sigma^2g(N) + \frac{2A}{R_0^2} }\right\}.
$$
 
Чтобы оценить зазор двойственности\ev{,} используем \eqref{eq:1}, также отметим, что данная оценка верна для любого $\boldsymbol{\lambda} \, \in \, \RR_+^m$. Поэтому\ev{,} беря минимум по \ev{всем $\boldsymbol{\lambda}$ из множества}  $\Lambda_{2R} =  \{\boldsymbol{\lambda} \, \in \, \mathbb{R}_{+}^m: \|\boldsymbol{\lambda}\|_2 \leq 2 R \}
$, получаем 
\begin{align*}
N \varphi(\boldsymbol{\hat{\lambda}}^N) \,  \le \,
\min_{\boldsymbol{\lambda} \, \in \, \Lambda_{2 R}}
\left \{
\sum_{t = 0}^{N-1} 
\left (
\varphi(\boldsymbol{\lambda}^t) +
\langle \nabla \varphi(\boldsymbol{\lambda}^t, \, \xi^t),
\boldsymbol{\lambda} - \boldsymbol{\lambda}^t \rangle
\right ) 
+  \frac{1}{2 \beta}
 \| \boldsymbol{\lambda}^0 - \boldsymbol{\lambda} \|_2^2
\right \}
\\
+ \, \frac{N\beta M^2}{2} , \,
\end{align*}
где для оценки последнего слагаемого использовалось   \ev{предположение}~\eqref{eq_stoch_gradM}.
Также учитывалось, что $\| \boldsymbol{\lambda}^{N} - \boldsymbol{\lambda} \|^2_2 \geq 0$.
В силу~\eqref{eq_fgm_ll0} получаем следующую оценку\ev{:}
$$
 \varphi(\boldsymbol{\hat{\lambda}}^N)  \, \le \,
 \frac{1}{N}
\min_{\boldsymbol{\lambda} \, \in \, \Lambda_{2R}}
\left \{
\sum_{t = 0}^{N-1} 
\left (
\varphi(\boldsymbol{\lambda}^t) +
\langle \nabla \varphi(\boldsymbol{\lambda}^t, \, \xi^t),
\boldsymbol{\lambda} - \boldsymbol{\lambda}^t \rangle
\right ) 
\right \}
+ 
\frac{2R^2}{ \beta N} +
\frac{\beta M^2}{2}.
$$
Прибавим и вычт\ev{е}м из выражения под минимумом $\sum\limits_{t = 0}^{N-1} \langle \nabla \varphi(\boldsymbol{\lambda}^t), \,
\boldsymbol{\lambda} - \boldsymbol{\lambda}^t \rangle$.
\ev{Тогда получаем} 
\begin{align*}
\min_{\boldsymbol{\lambda} \, \in \, \Lambda_{2 R}} &
\left \{
\sum_{t = 0}^{N-1} 
\left (
\varphi(\boldsymbol{\lambda}^t) +
\langle \nabla \varphi(\boldsymbol{\lambda}^t), \,
\boldsymbol{\lambda} - \boldsymbol{\lambda}^t \rangle
\right ) 
\right \}
\\
& \le \, 
 \min_{\boldsymbol{\lambda} \, \in \, \Lambda_{2 R}} 
\left \{
\sum_{t = 0}^{N-1} 
\left (
\varphi(\boldsymbol{\lambda}^t) +
\langle \nabla \varphi(\boldsymbol{\lambda}^t, \, \xi^t), \,
\boldsymbol{\lambda} - \boldsymbol{\lambda}^t \rangle
\right ) 
\right \} \\&  
+ 
\max_{\boldsymbol{\lambda} \, \in \, \Lambda_{2 R}}
\left \{
\sum_{t = 0}^{N-1} \langle \nabla \varphi(\boldsymbol{\lambda}^t, \, \xi^t)  -
\nabla \varphi(\boldsymbol{\lambda}^t), \,  \boldsymbol{\lambda}\rangle 
\right \}  
\\ &
+ \sum_{t = 0}^{N-1} \langle \nabla \varphi(\boldsymbol{\lambda}^t, \, \xi^t) -
\nabla \varphi(\boldsymbol{\lambda}^t), \, - \boldsymbol{\lambda}^t \rangle.
\end{align*}
Заметим, что $ - \boldsymbol{\lambda}^* \, \in \, \Lambda_{2 R}$. Тогда имеем
\begin{eqnarray*}
\sum_{t = 0}^{N-1} \langle \nabla \varphi(\boldsymbol{\lambda}^t, \, \xi^t) -
\nabla \varphi(\boldsymbol{\lambda}^t), \, - \boldsymbol{\lambda}^t \rangle & = &\sum_{t = 0}^{N-1}  \langle \nabla \varphi(\boldsymbol{\lambda}^t, \, \xi^t) -
\nabla \varphi(\boldsymbol{\lambda}^t), \, \boldsymbol{\lambda}^* - \boldsymbol{\lambda}^t  \rangle \\ 
& + & \sum_{t = 0}^{N-1}  \langle \nabla \varphi(\boldsymbol{\lambda}^t, \, \xi^t) -
\nabla \varphi(\boldsymbol{\lambda}^t), \, - \boldsymbol{\lambda}^* \rangle  \\ 
& \leq &
\max_{\boldsymbol{\lambda} \, \in \, \Lambda_{2 R}}
\sum_{t = 0}^{N-1}
 \langle \nabla \varphi(\boldsymbol{\lambda}^t, \, \xi^t)  -
\nabla \varphi(\boldsymbol{\lambda}^t), \,   \boldsymbol{\lambda}\rangle 
\\
&  + & \sum_{t = 0}^{N-1} \langle \nabla \varphi(\boldsymbol{\lambda}^t, \, \xi^t) - \nabla \varphi(\boldsymbol{\lambda}^t), \, \boldsymbol{\lambda}^* - \boldsymbol{\lambda}^t  \rangle.
\end{eqnarray*}
Отсюда получаем следующую оценку\ev{:} 
\begin{subequations}\label{ap_al_philamn}
\begin{align*}
\begin{split}
 \varphi(\boldsymbol{\hat{\lambda}}^N) &  \, \le \,
 \frac{1}{N}
\min_{\boldsymbol{\lambda} \, \in \, \Lambda_{2R}}
\left \{
\sum_{t = 0}^{N-1} 
\left (
\varphi(\boldsymbol{\lambda}^t) +
\langle \nabla \varphi(\boldsymbol{\lambda}^t, \, \xi^t), \,
\boldsymbol{\lambda} - \boldsymbol{\lambda}^t \rangle
\right ) 
\right \}
+ 
\frac{2 R^2}{ \beta N} +
\frac{\beta M^2}{2}  
\end{split}
\\
\begin{split}
& \leq 
\frac{1}{N}
\min_{\boldsymbol{\lambda} \, \in \, \Lambda_{2 R}}
\left \{
\sum_{t = 0}^{N-1} 
\left (
\varphi(\boldsymbol{\lambda}^t) +
\langle \nabla \varphi(\boldsymbol{\lambda}^t), \,
\boldsymbol{\lambda} - \boldsymbol{\lambda}^t \rangle
\right ) 
\right \}
\end{split} \\
\begin{split}
& + \frac{1}{N} \sum_{t = 0}^{N-1} \langle \nabla \varphi(\boldsymbol{\lambda}^t, \, \xi^t)-  \nabla \varphi(\boldsymbol{\lambda}^t), \, \boldsymbol{\lambda}^* - \boldsymbol{\lambda}^t \rangle \end{split}
\\
\begin{split}
& + 
\frac{2}{N} \max_{\boldsymbol{\lambda} \, \in \, \Lambda_{2 R}}
\left \{
\sum_{t = 0}^{N-1} \langle \nabla \varphi(\boldsymbol{\lambda}^t, \, \xi^t)  -
\nabla \varphi(\boldsymbol{\lambda}^t), \, \boldsymbol{\lambda}\rangle 
\right \}  + 
\frac{2 R^2}{ \beta N} +
\frac{\beta M^2}{2}. 
\end{split}
\tag{\ref{ap_al_philamn}}
\end{align*}
\end{subequations}
Из определения нормы получаем\ev{, что} 
$$
\max_{\boldsymbol{\lambda} \, \in \, \Lambda_{2 R}}
\left \{
\sum_{t = 0}^{N-1} \left \langle \nabla \varphi(\boldsymbol{\lambda}^t, \, \xi^t)  -
\nabla \varphi(\boldsymbol{\lambda}^t), \, \boldsymbol{\lambda} \right \rangle 
\right \} 
$$
$$
\leq 2 R \left \|
\sum_{t = 0}^{N-1} \left ( 
\nabla \varphi(\boldsymbol{\lambda}^t, \, \xi^t) -
\nabla \varphi(\boldsymbol{\lambda}^\ev{l}) 
\right ) 
\right \|_2.
$$
Используя \eqref{estimation_1}, получаем\ev{,} что с вероятностью $1-\delta$ выполняется 
\begin{equation}
\label{est_part_1}
   \left \|
\sum_{t = 0}^{N-1} 
\left ( 
\nabla \varphi(\boldsymbol{\lambda}^t, \, \xi^t) -
\nabla \varphi(\boldsymbol{\lambda}^\ev{l}) 
\right ) 
\right \|_2 \leq \sigma \sqrt{2 N} 
 \left (1 + \sqrt{3 \ln \frac{1}{\delta}}
 \right ).  
\end{equation}
\evv{Подставим в выражение $\sum\limits_{t = 0}^{N-1} \left (
\varphi(\boldsymbol{\lambda}^t) +
\langle \nabla \varphi(\boldsymbol{\lambda}^t), \,
\boldsymbol{\lambda} - \boldsymbol{\lambda}^t \rangle
\right )$ из
\eqref{ap_al_philamn}
значения $\varphi(\boldsymbol{\lambda}^t)$
 и $\nabla \varphi(\boldsymbol{\lambda}^t)$,
 получим
$$
\sum\limits_{t = 0}^{N-1} \left (
\langle \boldsymbol{\lambda^t}, \mathbf{b} \rangle + \sum_{k = 1}^n (u_k(x_k(\boldsymbol{\lambda^t})) - \langle \boldsymbol{\lambda}^t, \, \mathbf{C}_k x_k(\boldsymbol{\lambda^t}) \rangle ) +
\langle\mathbf{b} - C\mathbf{x}^t (\boldsymbol{\lambda^t}), \,
\boldsymbol{\lambda} - \boldsymbol{\lambda}^t \rangle
\right ) =
$$
$$
= \sum\limits_{t = 0}^{N-1} \left (
\sum_{k = 1}^n (u_k(x_k(\boldsymbol{\lambda^t})) +
\langle\mathbf{b} - C\mathbf{x}^t (\boldsymbol{\lambda^t}), \,
\boldsymbol{\lambda} \rangle
\right ).
$$
Тогда в силу вогнутости функций $u_k(x_k)$
$$
\frac{1}{N}
\min_{\boldsymbol{\lambda} \, \in \, \Lambda_{2 R}}
\left \{
\sum_{t = 0}^{N-1} 
\left (
\varphi(\boldsymbol{\lambda}^t) +
\langle \nabla \varphi(\boldsymbol{\lambda}^t), \,
\boldsymbol{\lambda} - \boldsymbol{\lambda}^t \rangle
\right ) 
\right \}
$$
$$
\leq U(\mathbf{\hat{x}}^N) 
- \frac{1}{N}
\max_{\boldsymbol{\lambda} \, \in \, \Lambda_{2 R}} 
\left \{
\sum_{t = 0}^{N-1} 
\langle  C \mathbf{x}^t(\boldsymbol{\lambda}^t) - \mathbf{b}, \,
\boldsymbol{\lambda} \rangle 
\right \}.
$$
Учитывая последнее неравенство, из \eqref{ap_al_philamn}
получаем} 
 \begin{align*}
 \varphi(\boldsymbol{\hat{\lambda}}^N) \, & \le \,
U(\evv{\mathbf{\hat{x}}^N})
- \frac{1}{N}
\max_{\boldsymbol{\lambda} \, \in \, \Lambda_{2 R}} 
\left \{
\sum_{t = 0}^{N-1} 
\langle  C \mathbf{x}^t(\boldsymbol{\lambda}^t) - \mathbf{b}, \,
\boldsymbol{\lambda} \rangle 
\right \} + \frac{2R^2}{ \beta N} +
\frac{\beta M^2}{2}
 \\ & 
+ \frac{2 R}{N} 
\left \|
\sum_{t = 0}^{N-1} 
\left ( \nabla \varphi(\boldsymbol{\lambda}^t, \, \xi^t) -
\nabla \varphi(\boldsymbol{\lambda}^t) 
\right ) 
\right \|_2 
\\ &
+ \frac{1}{N} \sum_{t = 0}^{N-1} \langle \nabla \varphi(\boldsymbol{\lambda}^t, \, \xi^t)-  \nabla \varphi(\boldsymbol{\lambda}^t), \, \boldsymbol{\lambda}^* - \boldsymbol{\lambda}^t \rangle . 
\end{align*}
Отсюда, учитывая оценк\ev{y} \eqref{est_part_1} и результат~\eqref{estimation_1}, получаем\ev{,} что с вероятностью $1-3\delta$

 \begin{eqnarray}
 \label{part_1}
 \varphi(\boldsymbol{\hat{\lambda}}^N) 
 - U(\mathbf{\hat{x}}^N)  + 
2 R
\left \|
\left [C \hat{\mathbf{x}}^N - \mathbf{b}
\right ]_+ 
\right \|_2  \, & \le & \,  \frac{2 R \sigma \sqrt{2} 
 \left (1 + \sqrt{3 \ln \frac{1}{\delta}}
 \right )}{\sqrt{N}}
 \notag \\
& + & \frac{2 R^2}{\beta N} +
\frac{\beta M^2}{2}   +    C_1\dfrac{\sigma\sqrt{ g(N)J} R ^ 2}{\sqrt{N}}. 
 \end{eqnarray}

По теореме~2.1 из~\cite{JudNem} для всех $\gamma > 0$ имеем
\begin{eqnarray*}
\mathbb{P} 
\left \{ 
\left \| \sum\limits_{t=0}^{N-1} \left( \mathbf{x}(\boldsymbol{\lambda}^t, \, \xi^t)-\mathbf{x}(\boldsymbol{\lambda}^t)\right) \right \|_2 \, \ge \,
(\sqrt{2} + \sqrt{2} \gamma) \sqrt{ \sum\limits_{t=0}^{N-1} \sigma_x^2} \, | \,
\{\xi^t\}_{t=0}^{N-1} 
\right \} 
\\
\, \le \,
\exp \left(-\frac{\gamma^2}{3} \right ).
\end{eqnarray*}
Выбирая $\gamma=\sqrt{3\ln \frac{1}{\delta}}$, получаем\ev{,} что с вероятностью $1 - \delta$
\begin{eqnarray*}
\|\Tilde{\mathbf{x}}^N - \hat{\mathbf{x}}^N \|_2  =  \dfrac{1}{N} \left \| \sum\limits_{t=0}^{N-1} \left( \mathbf{x}(\boldsymbol{\lambda}^t, \xi^t)-\mathbf{x}(\boldsymbol{\lambda}^t)\right) \right \|_2 
 \leq  \sigma_x \sqrt{\frac{2 }{N}} 
 \left (1 + \sqrt{3 \ln \frac{1}{\delta}}
 \right ).
\end{eqnarray*}
Тогда с вероятностью $1 - \delta$ выполняется следующее неравенство\ev{:} 
\begin{eqnarray*}
\|C\Tilde{\mathbf{x}}^N - C\hat{\mathbf{x}}^N \|_2  \leq \|C\|_2\cdot \|\Tilde{\mathbf{x}}^N - \hat{\mathbf{x}}^N \|_2 
 \leq  \sigma_x \sqrt{\frac{2 \lambda_{max}\left( C^{\mbox{T}}C \right)}{N}}
 \left (1 + \sqrt{3 \ln \frac{1}{\delta}}
 \right ).
\end{eqnarray*}
Заметим, что верно следующее:
\begin{eqnarray}
\label{part_2}
2 R
\left \|
\left [C \Tilde{\mathbf{x}}^N - \mathbf{b}
\right ]_+ 
\right \|_2 & =& \max_{\boldsymbol{\lambda} \, \in \, \Lambda_{2 R}} 
\left \{
\langle  C \Tilde{\mathbf{x}}^N - \mathbf{b}, \,
\boldsymbol{\lambda} \rangle 
 + \langle  C \hat{\mathbf{x}}^N -  C \hat{\mathbf{x}}^N - \mathbf{b} + \mathbf{b}, \,
\boldsymbol{\lambda} \rangle 
\right \} \notag\\
&\leq & 
\max_{\boldsymbol{\lambda} \, \in \, \Lambda_{2 R}} 
\left \{
\langle  C \hat{\mathbf{x}}^N - \mathbf{b}, \,
\boldsymbol{\lambda} \rangle 
\right \} + \max_{\boldsymbol{\lambda} \, \in \, \Lambda_{2 R}} 
\left \{
\langle  C \Tilde{\mathbf{x}}^N -  C \hat{\mathbf{x}}^N, \,
\boldsymbol{\lambda} \rangle 
\right \} \notag
\\
& \leq & 2 R
\left \|
\left [C \hat{\mathbf{x}}^N - \mathbf{b}
\right ]_+ 
\right \|_2 + 2 R\|C\Tilde{\mathbf{x}}^N - C\hat{\mathbf{x}}^N \|_2  \notag \\
& \leq & 2 R
\left \|
\left [C \hat{\mathbf{x}}^N - \mathbf{b}
\right ]_+ 
\right \|_2 
\notag \\
& + & 2 R\sigma_x \sqrt{\frac{2 \lambda_{max}\left( C^{\mbox{T}}C \right) }{N}} 
 \left (1 + \sqrt{3 \ln \frac{1}{\delta}}
 \right ).
\end{eqnarray}
В силу липшицевости функции $U$ получаем 
\begin{eqnarray*}
|U(\Tilde{\mathbf{x}}^N) - U(\hat{\mathbf{x}}^N) |  \leq M_{U} \|\Tilde{\mathbf{x}}^N - \hat{\mathbf{x}}^N \|_2 \leq M_{U} \sigma_x \sqrt{\frac{2}{N}} 
 \left (1 + \sqrt{3 \ln \frac{1}{\delta}}
 \right ).
\end{eqnarray*}
Тогда имеем
\begin{eqnarray}
\label{part_3}
U(\hat{\mathbf{x}}^N) = U(\Tilde{\mathbf{x}}^N) + \left(U(\hat{\mathbf{x}}^N) - U(\Tilde{\mathbf{x}}^N) \right) \geq U(\Tilde{\mathbf{x}}^N) - M_{U} \sigma_x \sqrt{\frac{2}{N}} 
 \left (1 + \sqrt{3 \ln \frac{1}{\delta}}
 \right ).
\end{eqnarray}
Подставляя~\eqref{part_2} и \eqref{part_3} в \eqref{part_1}, получаем\ev{,} что с вероятность\ev{ю} $1-4\delta$ выполняется 
 \begin{eqnarray*}
 \varphi(\boldsymbol{\hat{\lambda}}^N) 
 & - & U(\Tilde{\mathbf{x}}^N)  + 
2 R
\left \|
\left [C \Tilde{\mathbf{x}}^N - \mathbf{b}
\right ]_+ 
\right \|_2   \le    C_1\dfrac{\sigma\sqrt{ g(N)J} R ^ 2}{\sqrt{N}} 
+ \frac{2 R^2}{\beta N} +
\frac{\beta M^2}{2}  \notag \\ & + &    \frac{ \sqrt{2} \left (1 + \sqrt{3 \ln \frac{1}{\delta}} \right )}{\sqrt{N}}  \left( M_{U} \sigma_x + 2R\left(\sigma + \sigma_x \sqrt{\lambda_{max}\left( C^{\mbox{T}}C \right)} \right) \right). 
 \end{eqnarray*}

\subsection*{Доказательство теоремы~\ref{th_ellipsoid}}
	    
	    Так как $\|\nabla \varphi(\boldsymbol{\lambda})\|_2 \leq M$
	        \ev{ для любых } $\boldsymbol{\lambda} \in \Lambda_{2R}$ \ev{(см.~\eqref{eq_gradM})}, то справедлива следующая оценка:
	    $$
	    \sup_{\boldsymbol{\lambda}^1, \, \boldsymbol{\lambda}^2 \, \in \, \Lambda_{2R}} \langle \nabla \varphi(\boldsymbol{\lambda}^1), \, \boldsymbol{\lambda}^2 - \boldsymbol{\lambda}^1 \rangle \leq M \cdot 4R.
	    $$
Из теоремы 4.1 \cite{nemirovski2010accuracy} получаем:
	   $$
	    \max_{\boldsymbol{\lambda} \, \in \, \Lambda_{2R}} \sum_{t=1}^N \xi^t \langle \nabla \varphi(\boldsymbol{\lambda}^t), \, \boldsymbol{\lambda}^t - \boldsymbol{\lambda} \rangle \leq \varepsilon_N,
	    $$
где \ai{$\varepsilon_N= 32 \cdot 4 M R \exp \left\{ -\frac{N}{2 m (m+1)} \right\}$.}
Тогда имеем
\begin{equation*}
    \forall \, \boldsymbol{\lambda} \, \in \, \Lambda_{2R} \, \, \, \, \, \, \,\, \sum_{t \,  \in \, I_{N}} \xi^t \langle \nabla \varphi(\boldsymbol{\lambda}^t), \, \boldsymbol{\lambda}^t - \boldsymbol{\lambda} \rangle \leq \sum_{t=1}^N \xi^t \langle \nabla \varphi(\boldsymbol{\lambda}^t), \, \boldsymbol{\lambda}^t - \boldsymbol{\lambda} \rangle \leq \varepsilon_N.
\end{equation*}

Отсюда получаем, что верна следующая оценка:
	   $$
	    \sum_{t \, \in \, I_{N}} \xi^t \langle \mathbf{b} - C \mathbf{x}^t, \, \boldsymbol{\lambda}^t \rangle + \max_{\boldsymbol{\lambda} \, \in \, \Lambda_R} \left\langle -\sum_{t \, \in \, I_{N}} \xi^t (\mathbf{b} - C \mathbf{x}^t), \, \boldsymbol{\lambda} \right\rangle \leq \varepsilon_N,
	    $$
которую можно переписать в следующем виде:
\begin{equation}
\label{eq:est_ellips}
      \sum_{t \, \in \, I_{N}} \xi^t \langle \mathbf{b} - C \mathbf{x}^t, \, \boldsymbol{\lambda}^t \rangle \leq \varepsilon_N - 2R \left \|[C \hat{\mathbf{x}}^N - \mathbf{b}]_{+} \right \|_{2}.
\end{equation}
Далее, \ev{в силу} \eqref{eq_xi_argmax}, для каждого $\mathbf{x} \geq 0$ \ev{и $t \, \in \, I_{N}$} выполнено
	    $$
	    U(\mathbf{x}^t(\boldsymbol{\lambda}^t)) - \langle C \mathbf{x}^t(\boldsymbol{\lambda}^t) - \mathbf{b}, \, \boldsymbol{\lambda}^t \rangle \geq U(\mathbf{x}) - \langle C \mathbf{x} - \mathbf{b}, \, \boldsymbol{\lambda}^t \rangle.
	    $$
Умножая $t$-е неравенство на $\xi^t$, суммируя по всем индексам из $I_{N}$ и учитывая, что $\sum\limits_{t \, \in \, I_{N}} \xi^t U(\mathbf{x}^t) \leq U(\hat{\mathbf{x}}^{N})$ в силу вогнутости функций $u_k(x_k)$, $k = 1, \, \ldots, \, N$, получаем

	   $$
	    U(\mathbf{x}) - U(\hat{\mathbf{x}}^N) + \langle \mathbf{b} - C \mathbf{x}, \, \hat{\boldsymbol{\lambda}}^N \rangle \leq \sum_{t \in I_{N}} \xi^t \langle \mathbf{b} - C \mathbf{x}^t, \, \boldsymbol{\lambda}^t \rangle,
	    $$
 где $\hat{\boldsymbol{\lambda}}^N = \sum\limits_{t \, \in \, I_{N}} \xi^t \boldsymbol{\lambda}^t$. Используя оценку~\eqref{eq:est_ellips}, получаем
	    \begin{equation} \label{th_pr}
	    2R \left \|[C \hat{\mathbf{x}}^N - \mathbf{b}]_{+} \right \|_{2} + U(\mathbf{x}^{*}) - U(\hat{\mathbf{x}}^N) + \langle \mathbf{b} - C \mathbf{x}^{*}, \, \hat{\boldsymbol{\lambda}}^N \rangle \leq \varepsilon_N.
	    \end{equation}
Поскольку $\hat{\boldsymbol{\lambda}}^N \in \Lambda_{2R}$, и\ev{,} следовательно\ev{,} $\hat{\boldsymbol{\lambda}}^N \geq 0$, откуда $\langle \mathbf{b} - C \mathbf{x}^{*}, \hat{\boldsymbol{\lambda}}^N \rangle \geq 0$, 
\ev{из~\eqref{eq:est_ellips} следует}
$U(\mathbf{x}^{*}) - U(\hat{\mathbf{x}}^N) \leq \varepsilon_N$. 
Далее, так как для всех $\mathbf{x} \geq 0$, в силу определения $\boldsymbol{\lambda}^{*}$,
\ev{выполняется} $U(\mathbf{x}^{*}) \geq U(\mathbf{x}) - \langle \boldsymbol{\lambda}^{*}, \, C \mathbf{x} - \mathbf{b} \rangle$, получаем 
\begin{align*}
U(\hat{\mathbf{x}}^N)& \leq U(\mathbf{x}^{*}) - \langle \boldsymbol{\lambda}^{*}, \, \mathbf{b} - C \hat{\mathbf{x}}^N \rangle \leq  U(\mathbf{x}^{*}) - \min_{\boldsymbol{\lambda} \, \in
\, \lambda_{R}} \left\{\langle \boldsymbol{\lambda}, \, \mathbf{b} - C \hat{\mathbf{x}}^N \rangle \right\}=\\
&  U(\mathbf{x}^{*}) + \max_{\boldsymbol{\lambda} \, \in \, \lambda_{R}} \left\{\langle \boldsymbol{\lambda}, \, C \hat{\mathbf{x}}^N - \mathbf{b} \rangle \right\} \leq U(\mathbf{x}^{*}) + R \left \|[C \hat{\mathbf{x}}^N - \mathbf{b}]_{+} \right \|_{2}.
\end{align*}
Отсюда вместе с \eqref{th_pr} получаем $\evv{R} \|[C \hat{\mathbf{x}}^N - \mathbf{b}]_{+}\|_2 \leq \varepsilon_N$. Оценка~\eqref{th_saga} на количество итераций метода  следует из следующей выкладки:
$$
\varepsilon_N = 32 \cdot 4 M R \exp \left\{- \frac{N}{2m(m+1)}\right\} \leq \varepsilon \Rightarrow -\frac{N}{2m(m+1)} \leq \ln \left(\frac{\varepsilon}{32\cdot 4 M R}\right) 
$$
$$
\Rightarrow N \geq 2m(m+1) \ln \left(\frac{32 \cdot 4 M R}{\varepsilon}\right). 
$$


\begin{thebibliography}{99}
\bibitem{kelly}
\textit{Kelly F.P., Maulloo A.K., Tan D.K.H. Rate control for communication networks: shadow prices, proportional fairness and stability}~// Journal of the Operational Research Society. 1998. Vol.~49.  \textnumero~3. P.~237--252.

\bibitem{rokhlin}
\textit{Рохлин Д.Б. Распределение ресурсов в сетях связи с большим числом пользователей: стохастический метод градиентного спуска}~// Теория вероятностей и её применения (в печати).  2019.

\bibitem{arrow1958decentralization}
\textit{Arrow K.J., Hurwicz L.} Decentralization and computation in resource allocation. Stanford University, Department of Economics, 1958.

\bibitem{kakhbod2013resource}
\textit{Kakhbod A.} Resource allocation in decentralized systems with strategicagents: an implementation theory approach. Springer Science \& BusinessMedia, 2013.

\bibitem{campbell1987resource}
\textit{Campbell D.E.} Resource allocation mechanisms. Cambridge University Press, 1987.

\bibitem{friedman1995complexity}
\textit{Friedman E.J., Oren S.S. The complexity of resource allocation and price mechanisms under bounded rationality}~// Economic Theory.  1995. Vol.~6. \textnumero~2.  P.~225--250.

\bibitem{nesterov2018dual}
\textit{Nesterov Yu., Shikhman V. Dual subgradient method with averaging for optimal resource allocation}~// European Journal of Operational Research. 2018.  Vol.~270. \textnumero~3.  P.~907--916.


\bibitem{ivanova2018composite}
\textit{Ivanova A., Dvurechensky P., Gasnikov A., Kamzolov D.
Composite optimization for the resource allocation problem}~// arXiv preprint arXiv:1810.00595.  2018.

\bibitem{Nesterov1983}
\textit{Нестеров Ю.Е. Метод минимизации выпуклых функций со скоростью сходимости $O(1/k^2)$}~// Докл. АН СССР. 1983. Т.~269. \textnumero~39. C.~543--547.

\bibitem{gasnikov2016efficient}
\textit{Gasnikov A.V., Gasnikova E.V., Nesterov Yu.E.,  Chernov A.V.
Efficient numerical methods for entropy-linear programming problems}~// Comput. Math. and Math.Phys. 2016. Vol.~56.
\textnumero~4. P.~514--524.

\bibitem{chernov2016fast}
\textit{Chernov A., Dvurechensky P., Gasnikov A. Fast Primal-Dual Gradient Method for Strongly Convex Minimization Problems with Linear Constraints}~// Discrete Optimization and Operations Research: 9th International Conference, DOOR 2016, Vladivostok, Russia, September 19-23, 2016 Proceedings.  Springer International Publishing, 2016. P.~391--403.

\bibitem{dvurechensky2016primal-dual}
\textit{Dvurechensky P., Gasnikov A., Gasnikova E., Matsievsky S., Rodomanov A., Usik I.
Primal-Dual Method for Searching Equilibrium in Hierarchical Congestion Population Games}~// Supplementary Proceedings of the 9th International Conference on Discrete Optimization and Operations Research and Scientific School (DOOR 2016) Vladivostok, Russia, September 19-23, 2016. P.~584--595.  arXiv:1606.08988.

\bibitem{anikin2017dual}
\textit{Anikin A., Gasnikov A., Turin A., Chernov A.
Dual approaches to the minimization of strongly convex functionals with asimple structure under affine constraints}~// Comput. Math. and Math. Phys. 2017.  Vol.~57. \textnumero~8. P.~1262--1276.

\bibitem{dvurechensky2018computational}
\textit{Dvurechensky P., Gasnikov A., Kroshnin A. Computational Optimal Transport: Complexity by Accelerated Gradient Descent Is Better Than by Sinkhorn’s Algorithm}~// Proceedings of the 35th International Conference on Machine Learning. 2018.
Vol.~80. P.~1367--1376. arXiv:1802.04367.

\bibitem{nesterov2018primal-dual}
\textit{Nesterov Yu., Gasnikov A., Guminov S.,
Dvurechensky P.
Primal-dual accelerated gradient methods with small-dimensional relaxation oracle}~// arXiv:1809.05895. 2018. 

\bibitem{guminov2019acceleratedAM}
\textit{Guminov S., Dvurechensky P., Gasnikov A. On Accelerated Alternating Minimization}~// arXiv:1906.03622.  2019.

\bibitem{guminov2019accelerated}
\textit{Guminov S.V., Nesterov Yu.E., Dvurechensky P.E.,
Gasnikov A.V.
Accelerated Primal-Dual Gradient Descent with Linesearch for Convex, Nonconvex, and Nonsmooth Optimization Problems}~// Doklady Mathematics. 2019. Vol.~99.  P.~125--128.

\bibitem{kroshnin2019complexity}
\textit{Kroshnin A., Tupitsa N., Dvinskikh D., Dvurechensky P., Gasnikov A., Uribe C.A.
On the Complexity of Approximating Wasserstein Barycenters}~// Proceedings of the 36th International Conference on Machine Learning. 2019. Vol.~97. Eds. K.~Chaudhuri, R.~Salakhutdinov.  California, USA: PMLR, 2019. P.~3530--3540.  arXiv:1901.08686.

\bibitem{uribe2018distributed}
\textit{Uribe C.A., Dvinskikh D., Dvurechensky P., Gasnikov A., Nedich A.
Distributed Computation of Wasserstein Barycenters Over Networks}~// 2018 IEEE Conference on Decision and Control (CDC). 2018. P.~6544--6549.  arXiv:1803.02933.

\bibitem{dvinskikh2019primal}
\textit{Dvinskikh D., Gorbunov E., Gasnikov A.,
Dvurechensky P., Uribe C.A.
On Primal and Dual Approaches for Distributed Stochastic Convex Optimization over Networks}~// 2019 IEEE Conference on Decision and Control (CDC).  2019.  arXiv:1903.09844.

\bibitem{dvurechensky2018decentralize}
\textit{Dvurechensky P., Dvinskikh D.,
Gasnikov A., Uribe C.A., Nedic A.
Decentralize and Randomize: Faster Algorithm for Wasserstein Barycenters}~// Advances in Neural Information Processing Systems. 2018. Vol.~31. P.~10783--10793. arXiv:1806.03915.

\bibitem{Danskin}
\textit{Данскин Д.М.} Теория максмина. М.: Изд-во “Советское радио”, 1970.

\bibitem{Dem_book_74}
\textit{Демьянов В.Ф., Малоземов В. Н.} Введение в минимакс.  М.: Наука, 1972.

\bibitem{Nest}
\textit{Nesterov Yu. Smooth minimization of non-smooth functions}~// Mathematical Programming. 2005. Vol.~103. P.~127--152.

\bibitem{elip_nem}
\textit{Юдин Д.Б., Немировский А.С. Информационная сложность и эффективные методы решения выпуклых экстремальных задач}~// Экономика и мат. методы.  1976.  \textnumero~2. C.~357--369.

\bibitem{nemirovski2010accuracy}
\textit{Nemirovski A., Onn S., Rothblum U.G. Accuracy certificates for computational problems with convex structure}~// Mathematics of Operations Research. 2010. Vol.~35. \textnumero~1. P.~52--78.


\bibitem{lan2018random}
\textit{Lan G., Zhou Y. Random gradient extrapolation for distributed and stochastic optimization}~// SIAM Journal on Optimization.  2018. Vol.~28. \textnumero~4. P.~2753--2782.

\bibitem{bubeck}
\textit{Bubeck S. Convex optimization: Algorithms and complexity}~// Foundations and Trends in Machine Learning.  2015. Vol.~8. \textnumero~3/4. P.~231--357.

\bibitem{nesterov2018implementable}
Nesterov Yu. Implementable tensor methods in unconstrained convex optimization:  tech.  rep.  Universite  catholique  de  Louvain,  Center  for Operations Research and Econometrics (CORE). 2018.

\bibitem{gasnikov2019near}
\textit{Gasnikov A., Dvurechensky P., Gorbunov E.,  Vorontsova E., Selikhanovych D., Uribe C.A., Jiang B., Wang H., Zhang S., Bubeck S., Jiang Q., Lee Y.T., Li Y., Sidford A.
Near Optimal Methods for Minimizing Convex Functions with Lipschitz p-th Derivatives}~// Proceedings of the Thirty-Second Conference on Learning Theory. 2019. Vol.~99. P.~1392--1393.  URL: http://proceedings.mlr.press/v99/gasnikov19b.html. arXiv:1809.00382.

\bibitem{zhou2018simple}
\textit{Zhou K., Shang F., Cheng J. A simple stochastic variance reduced algorithm with fast convergence rates}~// arXiv preprint arXiv:1806.11027.  2018.

\bibitem{zhou2018direct}
\textit{Zhou K. Direct acceleration of SAGA using sampled negative momentum}~//arXiv preprint arXiv:1806.11048. 2018.

\bibitem{recht2011hogwild}
\textit{Niu F., Recht B., Re C., Wright S.J.
Hogwild: A lock-free approach to parallelizing stochastic gradient descent}~// Advances in Neural Information Processing Systems. 2011. P.~693--701.

\bibitem{jin2019short}
\textit{Jin C., Netrapalli P., Ge R., Kakade S.M., Jordan M.I. 
A short note on concentration inequalities for random vectors with subgaussian norm}~// arXiv preprint arXiv:1902.03736.  2019.

\bibitem{JudNem}
Juditsky A., Nemirovski A. Large  deviations  of  vector-valued  martin-gales  in  2-smooth  normed  spaces:   tech.  rep.   HAL: hal-00318071. 2008.
URL: http://hal.archives-ouvertes.fr/hal-00318071/. arXiv:0809.0813. 

\bibitem{nesterov_book2018}
\textit{Nesterov Yu.} Lectures on Convex Optimization. 2nd ed. Springer, 2018.

\bibitem{Nest_dis}
\textit{Нестеров Ю.Е.} Алгоритмическая выпуклая оптимизация: дисc.. д.ф.-м.н.: 01.01.07.  М.: МФТИ, 2013.

\bibitem{beck_matlab}
\textit{Beck A.} Introduction to Nonlinear Optimization:  Theory, Algorithms, and Applications with MATLAB. SIAM, Philadelphia, 2014.
\end{thebibliography}
\end{document}